\documentclass[12pt, twoside]{article}

\usepackage[french, english]{babel}
\usepackage{amsmath}
\usepackage{a4}
\usepackage{euscript}
\usepackage{amsfonts}
\usepackage{amssymb}
\usepackage{xypic}
\usepackage[dvips]{graphics}
\usepackage{enumerate}
\usepackage{array}
\usepackage{theorem}
\usepackage{epic}
\usepackage{color}  

\newtheorem{The}{Th\'eor\`eme}
\newtheorem{Pro}[The]{Proposition}
\newtheorem{Cor}[The]{Corollaire}
\newtheorem{Lem}[The]{Lemme}
\theorembodyfont{\upshape}
\newtheorem{Def}[The]{D\'efinition}
\newtheorem{Rem}[The]{Remarque}
\newtheorem{Ex}[The]{Exemple}

\def\Eproof{\hfill $\square$}

\newcommand{\Epl}{{\hfill $\lrcorner$ \medskip}}      
   
\def\NN{{\mathbb N}}
\def\ZZ{{\mathbb Z}}
\def\QQ{{\mathbb Q}}

\def\CC{{\mathbb C}}
\def\PP{{\mathbb P}}
\def\AA{{\mathbb A}}
\def\FF{{\mathbb F}}

\def\O{{\cal O}}

\def\C{{\cal C}}
\def\Y{{\cal Y}}
\def\M{{\cal M}}

\def\X{{\EuScript X}}
\def\Y{{\EuScript Y}}
\def\F{{\EuScript F}}

\def\p{{\mathfrak p}}
\def\q{{\mathfrak q}}

\def\GGa{\ensuremath{\mathbb G}_a}
\def\GGm{\ensuremath{\mathbb G}_m}
\def\mmu{\ensuremath{\boldsymbol \mu}}
\def\aalpha{\ensuremath{\boldsymbol \alpha}}

\def\gsigma{\langle\sigma\rangle}

\newcommand{\Spec}{{\mathrm{Spec}}}
\newcommand{\Spf}{{\mathrm{Spf}}}
\newcommand{\Ord}{{\mathrm{Ord}}}
\newcommand{\Pic}{{\mathrm{Pic}}}
\newcommand{\Res}{{\mathrm{Res}}}
\newcommand{\mult}{{\mathrm{mult}}}

\newcommand{\supp}{{\mathrm{supp}}}
\newcommand{\cartier}{{\mathrm{C}}}

\newcommand{\der}{{{d}}}
\newcommand{\Frac}{{\mathrm{Frac}}}
\newcommand{\height}{{\mathrm{ht}}}

\newcommand{\Aut}{{\mathrm{Aut}}}
\newcommand{\car}{{\mathrm{car}}}
\renewcommand{\H}{{\mathrm{H}}}

\begin{document}    

\title{Espace de modules des courbes hyperelliptiques stables et une in\'egalit\'e de Cornalba-Harris-Xiao}

\author{Sylvain Maugeais}
\date{}
\maketitle

\selectlanguage{french}
\begin{abstract}
Dans cet article, on d\'emontre une in\'egalit\'e entre les degr\'es du conducteur et du fibr\'e de Hodge d'une courbe stable hyperelliptique sur un anneau de valuation discr\`ete. Pour cela, on est amen\'e \`a d\'emontrer des th\'eor\`emes de rel\`evement de rev\^etements $p$-cycliques afin d'obtenir des informations sur le bord de l'espace des modules des courbes hyperelliptiques. Le r\'esultat d\'ecoule alors de la g\'eom\'etrie de ce dernier. 
\end{abstract}

\selectlanguage{english}
\begin{abstract}
In this paper, we show an inequality between the degrees of the conductor and of the Hodge line bundle of a hyperelliptic curve over a discrete valuation ring. For that we show a result about the lifting of $p$-cyclic covers of the projective line. We then apply this result to get information about the boundary of the moduli space of hyperelliptic stable curves and show the inequality.
\end{abstract}

\selectlanguage{french}
{\def\thefootnote{\relax}
\footnote{{\bf Math. Subj. Class. (2000)} : 14G20, 14H10, 14H30, 14D15.}
\footnote{{\bf Mots cl\'es } : Courbes hyperelliptiques, courbes stables, sch\'ema de modules, rev\^etement, formal patching.}

\tableofcontents

\section{Introduction}

Soit $\C \rightarrow S$ une famille de courbes stables de genre $g\ge 2$ 
sur une base $S$. Il lui correspond alors une section $\pi : S \rightarrow 
\overline{M}_{g,S}$ o\`u $\overline{M}_{g,S}$ est l'espace des modules des courbes stables de genre $g$ sur $S$. Notons $\lambda$ le fibr\'e de Hodge sur $\overline{M}_{g,S}$
et $\delta$ le diviseur du bord $\overline{M}_{g,S}\setminus {M}_{g,S}$. 
Lorsque $S$ est une courbe projective
lisse sur $\CC$, Cornalba, Harris \cite{Co-Ha} et Xiao \cite{Xia} ont 
d\'emontr\'e ind\'ependamment l'in\'egalit\'e 
\begin{equation}\label{chx}
(8g+4)\deg \pi^* \lambda  \ge  g \deg \pi^* \delta. 
\end{equation}
Notre but est d'\'etudier l'analogue de cette in\'egalit\'e lorsque 
$S$ est le spectre d'un anneau de valuation discr\`ete. Dans cette
situation, tout faisceau inversible sur $S$ est isomorphe \`a ${\cal O}_S$,
et le degr\'e n'est donc pas d\'efini. En revanche, un diviseur de 
Cartier sur $S$ correspond \`a un \'el\'ement de $f\in K(S)^*$ modulo 
${\cal O}_S^*$, on peut donc lui associer canoniquement un degr\'e qui
est la valuation de $f$. Le membre de gauche de l'in\'egalit\'e (\ref{chx})
a un sens si on peut trouver un diviseur de Cartier $D$ tel que 
${\cal O}(D)\simeq \lambda$, et dont le support ne contient pas
$\pi(S)$ de sorte que $\pi^*D$ soit d\'efini. Le membre de droite
est bien d\'efini si $\pi(S)$ n'est pas contenu dans $\delta$,
ce qui \'equivaut \`a dire que $\C\to S$ est \`a fibre g\'en\'erique lisse. 

Dans ce travail, nous consid\'erons les courbes stables ${\C}\to S$ 
\`a fibre g\'en\'erique lisse et hyperelliptique. Supposons pour un instant que la fibre g\'en\'erique $\C_\eta$ de ce morphisme est un rev\^etement de degr\'e 2 de 
$\PP^1_{K(S)}$. 
On a une section rationnelle 
canonique du fibr\'e $(\pi^*\lambda)^{\otimes 8g+4}$ sur $S$, \`a savoir
$$\Lambda=\Delta^g \left ( {dx \over 2y+Q} \wedge \ldots \wedge{x^{g-1}dx \over 2y+Q}\right )^{\otimes 8g+4} \in \left ( \bigwedge^g 
\H^0({\C_{\eta}}, \Omega_{\C_{\eta}/K}) \right ) ^{\otimes 8g+4}.$$
o\`u $y^2+Q(x)y-P(x)$ est une \'equation affine de $\C_{\eta}$,
et o\`u $\Delta$ est le discriminant de l'\'equation comme d\'efini par exemple dans \cite{Liu2}, section 1.2. Lorsque 
$2$ est inversible dans $S$, Kausz a montr\'e 
(cons\'equence imm\'ediate de \cite{Kau}, Theorem 3.1) que
$$g^2 \delta_s\ge \Ord_s(\Lambda) \ge g \delta_s,$$
o\`u $\delta_s$ est la somme des \'epaisseurs des points doubles de la fibre sp\'eciale (ou, de mani\`ere \'equivalente, le nombre de points 
singuliers de la fibre sp\'eciale de la d\'esingularisation minimale de 
$\C$ sur $S$, cf \cite{MB}, proposition 2.4). L'in\'egalit\'e de droite 
est l'analogue de (\ref{chx}) dans le cas local. 
La premi\`ere motivation 
de ce travail est le th\'eor\`eme suivant qui supprime l'hypoth\`ese 
$2$ inversible. 
\begin{The}
\label{ThmPrinc}
Soient $R$ un anneau de valuation discr\`ete, $S=\Spec(R)$, et  
$\C \rightarrow S$ une courbe stable dont la fibre g\'en\'erique $C$ est 
lisse et hyperelliptique. Notons $\Lambda$ la section canonique de 
$\left ( \bigwedge^g \H^0(C, \Omega_{C/K}) \right ) ^{\otimes 8g+4}$, 
on a alors
\begin{equation}\label{TheoPrinc}
g^2 \delta_s\ge \Ord_s(\Lambda)  \ge g\delta_s. 
\end{equation}
\end{The}

Notons que dans le cas d'une courbe elliptique \`a r\'eduction 
semi-stable, ces in\'egali\-t\'es sont des \'egalit\'es. 
Une telle in\'egalit\'e permet d'obtenir une borne sur $\omega_{\C/S}^2$ 
dans le cas arithm\'etique.

\begin{Cor}
Fixons une m\'etrique $C^\infty$ $\|\ \|_{Mod}$ sur $\omega_{\C_g/\bar H_g}$ ($\bar H_g$ est l'espace des modules des courbes stables munies d'un plongement tricanonique et $\C_g$ est la courbe stable universelle sur $\bar H_g$) et un $\varepsilon > 0$.
Soient $K$ un corps de nombres, $S=\Spec(\O_K)$, $\pi:\C \rightarrow S$ une courbe stable de genre $g \ge 2$ \`a fibre g\'en\'erique lisse hyperelliptique. Alors il existe une constante $c'=c'(\|\  \|_{Mod}, \varepsilon)$ ne d\'ependant que de $\| \ \|_{Mod}$ et de $\varepsilon$ telle qu'on ait l'in\'egalit\'e suivante 
$$(\omega_{\C/S}.\omega_{\C/S}) \le \left ( 3 {g^2+\varepsilon \over 2g+1}-1\right ) \left ( \sum_{{\mathfrak p} \in S} \delta_{\mathfrak p} log N_{\mathfrak p}+\sum_{\sigma:K \rightarrow \CC} \delta_\sigma(\C_{\sigma})\right )+c'[K:\QQ]$$ o\`u $N_{\mathfrak p}$ d\'esigne le cardinal du corps r\'esiduel de $S$ en $\mathfrak p$ et o\`u $(\omega_{\C/S}.\omega_{\C/S})$ (autointersection du fibr\'e dualisant au sens d'Arakelov) et $\delta_\sigma(\C_{\sigma})$ ($\delta$ invariant de $\C$) sont d\'efinis \`a l'aide de la m\'etrique $\| \ \|_{Mod}$ sur $\omega_{\C/S}$.
\end{Cor}

Cette borne est prouv\'ee dans \cite{Kau} sous l'hypoth\`ese que 
$\mathcal C$ a bonne r\'eduction au dessus de $2$. Elle 
est de la forme conjectur\'ee par Moret-Bailly dans \cite{MB} et on sait qu'une telle in\'egalit\'e, si elle \'etait montr\'ee dans le cas g\'en\'eral (et pas seulement dans le cas hyperelliptique), aurait des cons\'equences arithm\'etiques importantes (tel que la conjecture $abc$). 

Le th\'eor\`eme \ref{ThmPrinc} peut \^etre consid\'er\'e comme la positivit\'e
des diviseurs $g^2\delta - \Lambda$ et $\Lambda -g\delta$ sur un espace de 
modules des courbes stables hyperelliptiques. Comme la positivit\'e 
se v\'erifie aux points de codimension $1$, elle d\'ecoulera des in\'egalit\'es
(\ref{TheoPrinc}), d\'ej\`a prouv\'ees en \'egale caract\'eristique $0$ 
(cf. Appendice A), si l'espace de modules en question est suffisamment 
``sympathique''. Plus pr\'ecis\'ement, nous aurons besoin du r\'esultat 
d'irr\'eductibilit\'e qui suit. Pour tout sch\'ema $T$, notons $I_{g, T}$ l'adh\'erence 
sch\'ematique de l'espace de modules des courbes lisses hyperelliptiques 
$I_{g, T}$ dans le sch\'ema de Hilbert $\overline{H}_{g, T}$ des courbes
stables tricanoniques sur $T$. Cette construction ne donne pas \`a proprement 
parler un espace de modules, elle n'est pas m\^eme fonctorielle. Nous avons
cependant un r\'esultat de compatibilit\'e des fibres. 

\begin{The} \label{Thm3}
Soit $S$ le spectre d'un anneau de valuation discr\`ete, de point ferm\'e $s$.
Alors on a l'\'egalit\'e des ensembles 
$$\bar I_{g, S} \times_ S \Spec \ \kappa(s) = \bar I_{g, s}.$$ 
En particulier,  $\bar I_{g, S} \times_S \Spec\  \kappa(s)$ est irr\'eductible.
\end{The}

Dans un article \`a venir, nous construirons un espace de modules fin des courbes 
stables hyperelliptiques et \'etudierons ses diverses propri\'et\'es 
(structure locale, topologie, ...). 
Le th\'eor\`eme \ref{ThmPrinc} est prouv\'e dans la section 2, en
admettant le th\'eor\`eme \ref{Thm3}. En terme plus concr\`et, ce dernier
peut s'\'enoncer comme suit. 

\begin{The} \label{Thm4}
Soient $k$ un corps alg\'ebriquement clos et $C \rightarrow \Spec\, k$ 
une courbe stable hyperelliptique 
{\rm (d\'efinition \ref{CourbeStableHyperelliptique})}. 
Alors il existe un anneau de valuation dis\-cr\`e\-te $R$ 
d'\'egale caract\'eristique, de corps r\'esiduel $k$, 
et une courbe $\C \rightarrow \Spec \ R$ dont la fibre g\'en\'erique est lisse et hyperelliptique et dont la fibre sp\'eciale s'identifie \`a $C$.
\end{The}

La preuve de ce th\'eor\`eme constitue le coeur de cet article. Esquissons 
bri\`eve\-ment la d\'emarche. La partie non triviale est quand $\car\, k=2$. 
Par d\'efinition, $C$ est la r\'eduction d'une courbe lisse hyperelliptique
en caract\'eristique nulle. On sait alors que $C$ est un rev\^etement $2$-admissible 
d'une courbe de genre arithm\'etique $0$ (d\'efinition \ref{RevpAdmissible}). 
{\it Grosso modo}, cela veut dire qu'il existe un morphisme fini surjectif
$C\to X$ de degr\'e $2$ avec $p_a(X)=0$, que chaque composante irr\'eductible 
de $C$ poss\`ede un ouvert non vide muni d'une structure de torseur sous un des 
groupes alg\'ebriques $\mathbb Z/2\mathbb Z$, $\mmu_2$ ou $\aalpha_2$ sur $k$, et
qu'aux points d'intersection de $C$ on se donne des invariants num\'eriques
v\'erifiant une condition de compatibilit\'e.  
On montre que quitte \`a ajouter des droites projectives \`a $C$ et \`a modifier 
la structure des torseurs, on peut munir $C$ d'un graphe de Hurwitz {\it adapt\'e} 
(d\'efinition \ref{DefAdapt}). Ensuite, on d\'eforme la partie lisse 
de $C$ \`a l'aide des th\'eor\`emes g\'en\'eraux de rel\`evement de torseurs
(section 3.2), et on d\'eforme formellement les points 
d'intersection. On applique alors la technique de ``formal patching'' pour 
relever globalement $C$. La d\'efinition d'un graphe de Hurwitz adapt\'e
assure que les conditions locales et globales de recollement sont bien 
v\'erifi\'ees. 

Le th\'eor\`eme \ref{Thm4} est \'egalement valable pour les rev\^etements 
$p$-admissibles d'une courbe de genre arithm\'etique $0$. La technique 
employ\'ee ici s'inpire des travaux d'Henrio \cite{HenrioNP} et de Sa{\"\i}di 
\cite{Sai1}, \cite{Sai2}, \cite{Sai3} 
sur le rel\`evement en in\'egales caract\'eristiques, avec des
probl\`emes sp\'ecifiques \`a la situation d'\'egale caract\'eristique. 

La preuve du th\'eor\`eme \ref{Thm4} est r\'epartie dans les sections 3 -- 5.
La section 3 commence (apr\`es quelques rappels sur les torseurs sous des sch\'emas en groupes d'ordre $p$) par d\'efinir ce qu'est une courbe hyperelliptique stable sur un corps alg\'ebriquement clos. Plus g\'en\'eralement, on donne la d\'efini\-tion d'un rev\^etement $p$-admissible en caract\'eristique $p$.
Le section 4 est d\'evou\'ee \`a la d\'emonstration par des m\'ethodes de ``formal patching'' d'un th\'eor\`eme de rel\`evement de certains rev\`etements $p$-admissibles (ne faisant intervenir que certains types de torseurs et munis de plus d'un arbre de Hurwitz) en \'egale caract\'eristique.
La section 5 permet de voir que tout rev\^etement $p$-admissible peut en fait \^etre transform\'e en un rev\^etement $p$-admissible satisfaisant les conditions de la section 4. On montre ainsi que les fibres g\'eom\'etriques de l'espace des modules sont irr\'eductibles.

L'annexe A donne une preuve simple des in\'egalit\'es (\ref{TheoPrinc}) dans
le cas o\`u $\car\, k \ne 2$, et o\`u $\mathcal C_s$ ne poss\`ede qu'un point 
singulier (donc au plus deux composantes irr\'eductibles). 
Le calcul est similaire \`a celui de \cite{Kau}, mais est plus simple d\^ue
\`a notre hypoth\`ese. 
Dans l'annexe $B$, on montre le th\'eor\`eme \ref{ThmPrinc} par un calcul explicite
lorsque $\mathcal C_s$ est comme ci-dessus, mais avec $\car\, k=2$.
Cet appendice n'est pas utilis\'e dans le reste de l'article.

L'auteur remercie Yannick Henrio pour ses conseils pr\'ecieux sur les m\'ethodes utilis\'ees pour le rel\`evement des courbes, et Qing Liu pour l'attention qu'il a port\'ee \`a mon travail.
\medskip

\noindent {\bf Convention}
Dans tout cet article, $R$ sera un anneau de valuation discr\`ete de valuation $\nu$, de corps des fractions $K$ et de corps r\'esiduel $k$, $t$ sera une uniformisante de $R$. Le point ferm\'e de $\Spec(R)$ sera  not\'e $s$ et le point g\'en\'erique $\eta$. Certains graphes seront orient\'es. Si $a$ est une ar\^ete orient\'ee d'un graphe orient\'e, on notera $\bar a$ l'ar\^ete oppos\'ee. Les exposants $^o$ et $^t$ tiendront alors respectivement pour le sommet origine et le sommet terminal. Les courbes sur un sch\'ema $S$ seront des courbes \`a fibres g\'eom\'etriquement connexes.

\section{R\'eduction \`a l'\'etude de l'espace des modules des courbes hyperelliptiques.}

On cherche dans cette section \`a d\'efinir un sch\'ema proche de ce que serait l'espace des modules des courbes stables hyperelliptiques s'il existait, et \`a en donner quelques propri\'et\'es. Plus pr\'ecisement, nous allons consid\'erer l'adh\'erence du lieu hyperelliptique lisse dans l'espace des modules des courbes stables.

\subsection{Construction du sch\'ema $\bar I_g$}

On fixe un entier $g\ge 2$. Pour tout sch\'ema $S$, consid\'erons le foncteur
$$\begin{array}{llll}
 & (\textrm{Sch\'emas}/S)& \rightarrow & (\textrm{Ensembles}) \\
& T & \mapsto & \left\{\begin{array}{l}
\pi : \C \rightarrow T, \textrm{courbe stable de genre $g$}\\
        \textrm {avec isomorphisme } \PP(\pi_*\omega_{\C/T}^{\otimes 3}) \cong \PP^{5g-6}_T
\end{array}\right\}/ \textrm{iso}    
\end{array}$$
On sait (\cite{De-Mu}) que ce foncteur est repr\'esentable par un sch\'ema de Hilbert $\bar {H}_{g,S}$. De m\^eme, le foncteur des courbes lisses de genre $g$ munies de plongement tricanonique est repr\'esentable par un sous sch\'ema ouvert de $\bar H_{g, S}$ qu'on notera $H_{g, S}$. La d\'efinition de courbe hyperelliptique lisse sur une base quelconque se trouve dans \cite{Lo-Kl}.

\begin{Def}
Soit $\C \rightarrow S$ une courbe lisse et projective de genre $g$. Elle est dite {\it hyperelliptique} s'il existe une involution $\sigma \in \Aut_S(\C)$ telle que pour tout point g\'eom\'etrique $\bar s$ de $S$ le quotient $\C_{\bar s} / <\sigma>$ soit isomorphe \`a $\PP^1_{\Spec \ \kappa(\bar s)}$.
\end{Def}

\begin{Pro}
\label{ProHyperelliptiqueLisse}
Soient $S$ un sch\'ema et $\pi:\C \rightarrow S$ une courbe projective lisse de genre $g \ge 2$.
\begin{enumerate}[(1)]
\item Si $\C \rightarrow S$ est hyperelliptique, alors il existe un morphisme fid\`element plat $T \rightarrow S$ et un morphisme fini et fid\`element plat $\C \times_S T \rightarrow \PP^1_T$ de degr\'e $2$.
\item Supposons que  $S$ est r\'egulier, connexe et de dimension $1$, de point g\'en\'erique $\eta$. Alors $\C \rightarrow S$ est hyperelliptique si et seulement si $C_\eta \rightarrow \Spec \ \kappa(\eta)$ est hyperelliptique.
\item La courbe $\pi:\C \rightarrow S$ est hyperelliptique si et seulement si l'image du morphisme canonique $\C \rightarrow \PP(\pi_*\Omega^1_{\C/S})$ est une courbe de genre 0 et si sa formation commute au changement de base.
\end{enumerate}
\end{Pro}

\noindent{\it Preuve :} 
La propri\'et\'e (1) provient de \cite{Lo-Kl} Theorem 5.5, la propri\'et\'e (2) de \cite{Lo-Kl} Proposition 5.14, et la propri\'et\'e (3) de \cite{Lo-Kl} Theorem 5.5. \Eproof
\medskip 

On souhaite alors prouver la repr\'esentabilit\'e du foncteur
$$\begin{array}{llll}
\underline{Hyp} :  & (\textrm{Sch\'emas}/S) & \rightarrow & (\textrm{Ensembles}) \\
& T & \mapsto & \left\{\begin{array}{l}
\pi : \C \rightarrow T, \textrm{courbe lisse}\\
	\textrm{hyperelliptique de genre $g$}\\
        \textrm{avec isomorphisme } \PP(\pi_*\omega_{\C/T}^{\otimes 3}) \cong \PP^{5g-6}_T
\end{array}\right\}/ \textrm{iso}    
\end{array}$$

\begin{Pro}
Le foncteur $\underline{Hyp}$ est repr\'esentable par un sous sch\'ema ferm\'e de 
${H}_{g,S}$.
\end{Pro}

\noindent{\it Preuve :}  On suit la d\'emonstration donn\'ee dans \cite{La-Lo}, Theorem 3. Soit $\pi : \mathfrak{C} \rightarrow {H_{g, S}}$ la courbe stable universelle de genre $g$  sur $S$, et $D$ l'image sch\'ematique du morphisme $f:\mathfrak{C} \rightarrow \PP(\pi_* \Omega^1_{\mathfrak{C}/{H_{g, S}}})=:\PP$ (on sait que $f$ est bien d\'efini car $\pi^* \pi_* \Omega^1_{\mathfrak{C}/{H_{g, S}}} \rightarrow \Omega^1_{\mathfrak{C}/{H_{g, S}}}$ est surjectif). On d\'efinit alors $Q$ comme \'etant le conoyau de l'homomorphisme induit $\O_\PP \rightarrow f_* \O_\mathfrak{C}$. Pour tout point g\'eom\'etrique $P$ de ${H_{g, S}}$, il ne peut se produire que deux cas :
\begin{enumerate}[(1)]
\item $f_P:{\mathfrak C}_P \rightarrow D_P$ est un isomorphisme et $Q_P=0$ (c'est \`a dire que le morphisme canonique $f$ est une immersion), auquel cas la courbe est non hyperelliptique;
\item $D_P \cong \PP^1_P$ et $Q_P \cong \O_{\PP^1}(-g-1)$ (c'est \`a dire que le morphisme canonique $f$ \`a pour image une droite projective), auquel cas la courbe est hyperelliptique
\end{enumerate}
(cf \cite{Lo-Kl}, section 5). Ainsi, seuls deux polyn\^omes de Hilbert peuvent arriver pour $Q_P$, donc la stratification platifiante de $Q$ sur ${H_{g, S}}$ (qui existe par \cite{Mum}, Lecture 8) est constitu\'ee de deux sous-sch\'emas disjoints localement ferm\'es $F_1$ et $F_2$ (correspondant respectivement aux cas $(1)$ et $(2)$). En particulier, $F_2$ est ferm\'e car il correspond au support de $Q$. Il reste \`a montrer que $F_2$ correspond au courbes hyperelliptiques mais ceci d\'ecoule de \cite{Lo-Kl} Proposition 2.13 et de la proposition \ref{ProHyperelliptiqueLisse}, 3  \Eproof    
\medskip

Pour tout sch\'ema $S$, on note alors ${I}_{g, S}$ le sous-sch\'ema ferm\'e de $H_{g, S}$ qui repr\'esente le foncteur $\underline{Hyp}$ (c'est \`a dire $F_2$). On notera $\bar{I}_{g, S}$ l'adh\'erence sch\'ema\-tique de ${I}_{g, S}$ dans $\bar {H}_{g, S}$. Si $S$ est r\'eduit (resp. int\`egre) alors $I_{g, S}$ aussi (car $I_{g, S} \rightarrow S$ est lisse par \cite{La-Lo} Theorem 2), et il en est donc de m\^eme pour $\bar I_{g, S}$. De plus, pour tout morphisme de sch\'ema $T \rightarrow S$, on a une immersion ferm\'ee $\bar I_{g, S} \times_S T \rightarrow \bar I_{g, T}$ et les deux sch\'emas ont m\^eme espaces sous jacents. Si $T \rightarrow S$ est plat alors ce morphisme est un isomorphisme, ceci d\'ecoulant du lemme

\begin{Lem}
Soit $S$ un sch\'ema, $f:X \rightarrow Y$ un morphisme de $S$ sch\'ema s\'epar\'e et de type fini. Notons ${\mathrm {Im}} f$ l'image sch\'ematique de $f$. Alors, pour tout sch\'ema morphisme plat $T \rightarrow S$ on a ${\mathrm {Im}} (f \times Id_T)={\mathrm {Im}} f \times_S T$.
\end{Lem}

\noindent{\it Preuve :} Le sch\'ema ${\mathrm {Im}} f$ est le sous sch\'ema ferm\'e de $Y$ d\'efinit par le faisceau d'id\'eaux $\ker \left (\O_Y \rightarrow f_* \O_X \right )$. Le r\'esultat provient alors du fait que la formation de $f_* \O_X$ commute au changement de base plat lorsque $f$ est s\'epar\'e et de type fini et, par suite, que la formation du noyau commute au changement de base plat. \Eproof

\begin{Rem}
\label{RepReduit}
Soit $T$ un $S$-sch\'ema int\`egre. Soit $\C \rightarrow T$ une courbe stable dont la fibre g\'en\'erique est lisse et hyperelliptique. On a alors un morphisme $T \rightarrow \bar H_{g, S}$ dont l'image est (ensemblistement) contenue dans $\bar I_{g, S}$. Par suite, ce morphisme se factorise par l'inclusion $\bar I_{g, S} \rightarrow \bar H_{g, S}$ car $T$ est r\'eduit.
\end{Rem}
 
Examinons maintenant le bord de $\bar I_{g, \CC}$. Soit $i$ un entier tel que $0 \le i \le g/2$. Soit $C \rightarrow \Spec \ k$ une courbe stable de genre $g$ sur un corps alg\'ebriquement clos $k$. Un point $P$ de $C$ sera dit de type $\beta_i$ si la normalisation de $C$ en $P$ est compos\'ee de deux composantes connexes de genre respectif $i$ et $g-i$ si $i \not = 0$, et connexe si $i=0$. Notons $\Delta_i$ le lieu des points de $\bar H_{g, \CC}$ correspondant \`a des courbes stables $C$ de genre $g$ ayant un point singulier de type $\beta_i$. Il est prouv\'e dans \cite{De-Mu}, \S 3 que  $\Delta_i$ est irr\'eductible et que $\bar H_{g, \CC} \setminus H_{g, \CC} = \cup_i \Delta_i$. D'autre part, on montre (\cite{Co-Ha}, Section 4, \S b) que si $i > 0$ alors $\Delta_i \cap \bar I_{g, \CC}$ est irr\'eductible, et qu'on a une d\'ecomposition en composantes irr\'eductibles  $$\Delta_0 \cap \bar I_{g, \CC}= \bigcup_{j=0}^{[(g-1)/2]} \Xi_j$$ o\`u $\Xi_j$ est d\'efinie de la mani\`ere suivante.
Soient $C \rightarrow \Spec \ \CC$ un point de $\bar I_{g, \CC}$ et $\sigma$ l'involution canonique de $C$ (qui provient du fait que $C$ est la r\'eduction d'une courbe hyperelliptique lisse). Si $j > 0$, un point singulier $P$ de $C$ est dit de type $\alpha_j$ si $P$ n'est pas fixe par $\sigma$ et que la normalisation partielle de $C$ en $P$ et $\sigma(P)$ est compos\'e de deux composantes connexes de genre respectif $j$ et $g-j-1$. Un point singulier $P$ est dit de type $\alpha_0$ s'il est fixe par $\sigma$ et si l'involution hyperelliptique permutte les branches analytiques de $C$ en $P$. L'ensemble $\Xi_j$ est alors le lieu des points de $\bar I_{g, \CC}$ repr\'esentant des courbes ayant un point de type $\alpha_j$.

\subsection{D\'emonstration des in\'egalit\'es (\ref{TheoPrinc})}

Soient $S$ un sch\'ema, $i:\bar I_{g, S} \rightarrow \bar H_{g, S}$ l'inclusion, $\delta$ le diviseur du bord (c'est \`a dire le diviseur obtenu par changement de base $S \rightarrow \Spec \ \ZZ$ \`a partir de $\bar H_{g, \ZZ} \setminus H_{g, \ZZ}$), on sait (\cite{De-Mu}, Corollary 1.9) que c'est un diviseur de Cartier, et $\lambda$ le fibr\'e de Hodge sur $\bar H_{g, S}$ (c'est \`a dire que pour une courbe stable $\pi:\C \rightarrow T$ \`a laquelle correspond un morphisme $\mu:T \rightarrow \bar H_{g, S}$ on a $\mu^* \lambda = \det \left( \pi_* \omega_{\C/T}\right)$). Dans le cas o\`u le sch\'ema de base est le spectre d'un anneau de valuation discr\`ete, on a les interpr\'etations suivantes :

\begin{Pro}
\label{IntDelta}
Soient $R$ un anneau de valuation discr\`ete dont le corps r\'esiduel est alg\'ebriquement clos et $S=\Spec \ R$. Soit $\C \rightarrow S$ un courbe stable dont la fibre g\'en\'erique est lisse. Soit $\mu:S \rightarrow \bar H_{g, S}$ une section correspondant \`a un plongement tricanonique de $\C$. Soient $q_1, \ldots, q_r$ les points singuliers de la fibre sp\'eciale $\C_s$ d'\'epaisseurs respectives $e(q_1), \ldots, e(q_r)$ dans $\C$, $t$ une \'equation locale de $\delta$ en $\mu(s)$. Alors $\sum_{i=1}^r e(q_i) = \nu_K(\mu^\#(t))$. 
\end{Pro}

\noindent{\it Preuve :}  Voir \cite{Liu3}, section 5, Lemme. \Eproof

\begin{Rem}
 Cette interpr\'etation est encore v\'erifi\'ee par l'image r\'ecipro\-que $\delta':=i^*\delta$ de $\delta$ sur $\bar I_{g, S}$.
\end{Rem}

\begin{Rem}
Notons $\lambda'$ l'image r\'eciproque du fibr\'e de Hodge $\lambda$ sur $\bar I_{g, S}$. Soient $\C \rightarrow S$ une courbe hyperelliptique et $T \rightarrow S$ un morphisme. On montre ais\'ement que la formation du fibr\'e $\lambda'$ commute au changement de base. 
\end{Rem}

\begin{Pro}
Soient $T$ un sch\'ema et $\C \rightarrow T$ une courbe lisse hyperelliptique. Alors il existe une section canonique $\Lambda_{\C/T} \in (\det \H^0(\C, \omega_{\C/T}))^{\otimes 8g+4}$ qui est une base et dont la formation commute au changement de base.
\end{Pro}

\noindent{\it Preuve :} Traitons d'abord le cas o\`u il existe un morphisme $f:\C \rightarrow \PP^1_T$ fid\`element plat et fini de degr\'e 2. Alors il existe un recouvrement de $T$ par des sch\'emas affines $T_i=\Spec \ A_i$ tel que $\C_{T_i}$ contient un ouvert isomorphe \`a 
$$\Spec \ A_i[x, y]/(y^2+Q_i(x)y-P_i(x)).$$
On montre alors que 
$$\Lambda_{\C_{T_i}/T_i}=\Delta^g \left ( {dx \over 2y+Q_i} \wedge \ldots \wedge{x^{g-1}dx \over 2y+Q_i}\right )^{\otimes 8g+4}$$
(o\`u $\Delta$ est le discriminant de l'\'equation comme d\'efini par exemple dans \cite{Liu2}, section 1.2) d\'efinit une base et ne d\'epend pas du choix de l'\'equation (cf par exemple \cite{Kau} Proposition 2.2).

Soit $V$ un ouvert affine de $T_i \cap T_j$, alors 
$$\Lambda_{\C_{T_i}/T_i}|V=\Lambda_{\C_{T_j}/T_j}|V$$ car ces deux sections ne d\'epende pas des \'equations choisies. On peut donc recoller les $\Lambda_{\C_{T_i}/T_i}$ et ainsi obtenir une section $\Lambda_{\C_{T}/T}$.

Dans le cas g\'en\'eral, d'apr\`es \cite{Lo-Kl}, Theorem 5.5, on peut trouver une famille de sch\'emas $T'_i$ affines et un morphisme fid\`element plat $\cup_i T_i \rightarrow T$ tel que pour tout $i$, il existe un morphisme fid\`element plat $\C \times_T T_i \rightarrow \PP^1_{T_i}$ de degr\'e 2. Par suite, on peut d\'efinir pour tout $i$ une section $\Lambda_{\C_{T'_i}/T'_i}$. Il reste alors \`a recoller ces sections. Mais comme la formation de $\lambda$ commute au changement de base, et que le pr\'efaisceau (fppf) associ\'e au faisceau (Zariski) coh\'erent $\lambda$ est un faisceau, le r\'esultat provient de l'unicit\'e de la section canonique sur $T_i \times_T T_j$. \Eproof
\medskip

Soient $S=\Spec \ R$ le spectre d'un anneau de valuation discr\`ete,  $\pi:\C \rightarrow S$ une courbe stable dont la fibre g\'en\'erique est lisse et hyperelliptique, et $\mu:S \rightarrow \bar I_{g, S}$ une section de $\bar I_{g, S}$ correspondante. Notons $\Lambda$ la section canonique $\Lambda_{\C_\eta /K}$ du fibr\'e de Hodge sur $S$. Si $\Lambda_0$ est une base de $\mu^* {\lambda'}^{\otimes 8g+4}$, il existe $a \in K$ tel que $\Lambda=a\Lambda_0$. On note alors $\Ord_s(\Lambda)=\nu(a)$ et $\delta_s$ la somme des \'epaisseurs des points doubles de la fibre sp\'eciale de $\C$.

Notons $\Lambda'$ la section canonique de ${\lambda'}^{\otimes 8g+4}$ sur $\bar I_{g, S}$ et $D$ le diviseur de Cartier sur $\bar I_{g, S}$ associ\'e \`a la section rationnelle $\Lambda'$.

\begin{Pro}
\label{IntD}
Soit $\mu:S \rightarrow \bar I_{g, S}$ un morphisme correspondant \`a une courbe stable hyperelliptique $\C \rightarrow S$ dont la fibre g\'en\'erique est lisse, alors $\mult (\mu^* D)=\Ord_s(\Lambda)$.
\end{Pro}

\noindent{\it Preuve :} Soit $\Lambda_0$ une base locale de ${\lambda'}^{\otimes 8g+4}$ au voisinage de $\mu(s)$. La fonction m\'eromorphe $a=\Lambda'/\Lambda_0$ est alors une \'equation locale du diviseur de Cartier $D$. On a, par d\'efinition, $\mult (\mu^* D)=\mult_R(\mu^*a)=\Ord_s(\Lambda)$. \Eproof

\begin{Lem}
\label{PosNorm}
Soit $A$ un anneau noetherien normal, alors pour tout $f \in \Frac(A)$, $f\in A$ si et seulement si pour tout $\p \in \Spec \ A$ de hauteur $1$ on a $\mult_{A_\p} (f) \ge 0$.
\end{Lem}

\noindent{\it Preuve :} Cela r\'esulte de l'\'egalit\'e $$A=\bigcap_{\p \in \Spec A, \height(\p)=1} A_{\p}.$$ \Eproof
\medskip 

Le th\'eor\`eme suivant, dont la d\'emonstration fait l'objet des sections suivantes, est alors le point clef de la d\'emonstration du 
th\'eor\`eme \ref{ThmPrinc}.

\begin{The}
\label{IgDense}
Soient $T$ un sch\'ema r\'egulier et de dimension 1, et $\tau$ un point de $T$. On a l'\'egalit\'e (ensembliste) $\bar I_{g, T} \times_T \Spec \ \kappa(\tau) = \bar I_{g, \Spec \kappa(\tau)}.$ En particulier,  $\bar I_{g, T}\times_T \Spec \ \kappa(\tau)$ est irr\'eductible.
\end{The}

Ce th\'eor\`eme sera prouv\'e dans la section 5.3.

\begin{Cor}
\label{WeilEff}
Soit $T$ un sch\'ema r\'egulier de dimension 1 dont le corps des fonctions est de caract\'eristique z\'ero ($T$ est donc excellent), alors
le diviseur de Cartier $D-g \delta '$ est effectif sur le normalis\'e ${\bar I}_{g, T}'$ de ${\bar I}_{g, T}$.
\end{Cor}

\noindent{\it Preuve :} D'apr\`es le lemme \ref{PosNorm}  il suffit de montrer que pour tout point g\'en\'erique $\xi'$ de ${\bar I}_{g, T}'$ de hauteur 1, on a $\mult_{\xi'}(D-g\delta') \ge 0$. Notons $\xi$ l'image de $\xi'$ dans $\bar I_{g, T}$. Comme $T$ est excellent, le morphisme ${\bar I}_{g, T}' \rightarrow \bar I_{g, T}$ est fini et $\bar I_{g, T}$ est universellement cat\'enaire. Par suite, d'apr\`es \cite{EGA} Proposition 5.6.10, le point $\xi$ est de hauteur 1 dans $\bar I_{g, T}$.

Notons $A$ l'anneau local de $\bar I'_{g, T}$ en $\xi'$ et $\C \rightarrow \Spec A$ la courbe correspondant au morphisme $\Spec \ A \rightarrow \bar I_{g, T}$. Ainsi d\'efini, $A$ est un anneau de valuation discr\`ete et la fibre g\'en\'erique de $\C \rightarrow \Spec A$ est lisse car le point g\'en\'erique de $A$ correspond au point g\'en\'erique de $\bar I_{g, T}$. Par suite, on a $\mult_{\xi'}(D-g\delta')=\Ord_{\xi'}(\Lambda)-g\delta_{\xi'}$.

 Supposons que $\xi'$ ne soit pas dans la fibre g\'en\'erique de ${\bar I}_{g, T}' \rightarrow T$. Notons $\tau$ l'image de $\xi'$ dans $T$. Comme $I_{g, T}\times_T \ \Spec \kappa(\tau)$ est dense dans $\bar I_{g, T}\times_T \Spec \ \kappa(\tau)$ d'apr\`es le th\'eor\`eme \ref{IgDense}, les points de $\bar I_{g, T}\times_T \Spec \ \kappa(\tau) \setminus I_{g, T}\times_T \Spec \ \kappa(\tau)$ sont de codimension $\ge 1$ dans $\bar I_{g, T}\times_T \Spec \ \kappa(\tau)$. Il s'ensuit que $\xi \in I_{g, \Spec \ \kappa(\tau)}$. Par suite, $\C$ est lisse et on sait que pour de telle courbe, on a $\Ord_{\xi'}(\Lambda)-g\delta_{\xi'}= 0$, c'est \`a dire que $D-g\delta'$ est effectif en $\xi'$. 

 Supposons que $\xi'$ soit dans la fibre g\'en\'erique de ${\bar I}_{g, T}' \rightarrow T$. Comme $\xi$ est de hauteur 1, la fibre sp\'eciale de $\C \rightarrow \Spec A$ poss\`ede au plus deux composantes irr\'eductibles (cf. fin de la section 2.1). D'apr\`es le calcul fait en caract\'eristique z\'ero (cf annexe A) on a $\Ord_s(\Lambda)-g \delta_s \ge 0$ donc $\mult_{\xi'}(D-g \delta') \ge 0$, ce qui prouve que $D-g\delta'$ est effectif sur ${\bar I}_{g, T}'$
 \Eproof

\begin{Lem}
\label{FoncReguliere}
Soient $X$ un sch\'ema int\`egre quasi-compact, $D_0$ un diviseur de Cartier effectif, $Z=\supp(D_0)$ et $U=X\setminus Z$. Soit $E$ un diviseur de Cartier sur $X$ tel que $E|_U$ soit effectif. Alors il existe un entier $n \in \NN$ tel que $D+nD_0$ soit effectif.
\end{Lem}

\noindent{\it Preuve :} Quitte \`a remplacer $X$ par des ouverts affines assez petits, on peut supposer que $X=\Spec \ A$. Soit $f$ une \'equation locale d\'efinissant $D_0$ et $s$ une \'equation locale d\'efinissant $D$. Comme on a $\O_X(U)=A_f$, $s$ est de la forme $a/f^n$, ce qui d\'emontre le lemme. \Eproof
\medskip

Soient $S$ le spectre d'un anneau de valuation discr\`ete et $\C \rightarrow S$ une courbe stable \`a fibre g\'en\'erique lisse hyperelliptique. Elle est alors munie d'un plongement tricanonique. Il lui correspond un morphisme $$\mu:S \rightarrow \bar H_{g, S}.$$ Comme la courbe $\C$ est hyperelliptique et que $R$ est r\'eduit, ce morphisme se factorise par $\bar I_{g, S}$ (cf Remarque \ref{RepReduit}) et on a le diagramme
$$\xymatrix{
S \ar[r]^{\mu}\ar[dr]^{\mu'} & \bar H_{g, S} \\
& \bar I_{g, S} \ar[u].}$$

\noindent{\it Preuve du th\'eor\`eme \ref{ThmPrinc}:} Tout d'abord, quitte \`a changer $R$ en son compl\'et\'e (ce qui ne change pas les valuations $\Ord_s(\Lambda)$ et $\delta_s$), on peut supposer que $R$ est complet. Par suite, il existe un sch\'ema excellent $T$ dont le corps de fraction est de carat\'eristique z\'ero et tel que la courbe $\C \rightarrow S$ corresponde \`a un morphisme $\mu:S \rightarrow \bar I_{g, T}$ (on peut prendre $T=\Spec R$ si $\car \Frac(R)=0$ et dans le cas contraire, $T=\Spec W(k)[[t]]$ o\`u $W(k)$ est l'anneau des vecteurs de Witt du corps r\'esiduel de $T$). Notons $x$ l'image du point ferm\'e de $S$ par $\mu$ et $A$ l'anneau local de $x$ dans $\bar I_{g, T}$. On a alors un morphisme $\mu^\# : A \rightarrow R$. Par suite, si $f$ est une \'equation locale de $D$ sur un voisinage $X=\Spec \ A$ de $x$, et $h$ est une \'equation locale de $g\delta'$, $f/h \in A'$ (cl\^oture int\'egrale de $A$) d'apr\`es le corollaire \ref{WeilEff} et le lemme \ref{PosNorm}. C'est \`a dire qu'il existe $a_{m-1}, \ldots, a_0 \in A$ tel que $(f/h)^m+a_{m-1}(f/h)^{m-1}+\ldots+a_0=0$. Or d'apr\`es le lemme \ref{FoncReguliere} il existe un entier $n$ tel que le diviseur de Cartier $D+n\delta'$ soit effectif, ce qui equivaut \`a dire que $fh^n$ est une fonction r\'eguli\`ere sur $A$. D'autre part, $\mu^\#(h) \not = 0$ (car la fibre g\'en\'erique de $\C$ est lisse) et comme $\mu^\#\left ( (fh^n)^m+a_{m-1}(fh^n)^{m-1}h^{n+1}+\ldots+a_0h^{(n+1)m} \right )=0$, on voit que $\mu^\#(fh^n)/\mu^\#(h^{n+1})$ est enti\`ere sur $R$, donc appartient \`a $R$. 

D'autre part, comme la formation de $D$ et de $\delta'$ commute au changement de base $S \rightarrow \Spec \ A$, l'ordre de $\mu^\#(fh^{n+1}/\mu^\#(h^n)$ dans $R$ est \'egal \`a $\Ord_s(\Lambda) - g \delta'$ d'apr\`es les propositions \ref{IntD} et \ref{IntDelta}, ce qui d\'emontre le th\'eor\`eme. 

Le m\^eme raisonnement pouvant \^etre tenu pour le diviseur $g^2\delta'-D$, 
les in\'egalit\'es
(\ref{TheoPrinc}) 
sont d\'emontr\'ees. \Eproof
\medskip

On peut d\`es lors supprimer l'hypoth\`ese de bonne r\'eduction en $2$ dans \cite{Kau} Theorem 6 et Theorem 7.1, et obtenir le corollaire suivant 
\begin{Cor}
Fixons une m\'etrique $C^\infty$ $\|\ \|_{Mod}$ sur $\omega_{\C_g/\bar H_g}$ ($\C_g$ est la courbe stable universelle sur $\bar H_g$) et un $\varepsilon > 0$.
Soient $K$ un corps de nombres, $S=\Spec(\O_K)$, $\pi:\C \rightarrow S$ une courbe stable de genre $g \ge 2$ dont la fibre g\'en\'erique est hyperelliptique. Alors il existe une constante $c'=c'(\|\  \|_{Mod}, \varepsilon)$ ne d\'ependant que de $\| \ \|_{Mod}$ et de $\varepsilon$ telle qu'on ait l'in\'egalit\'e suivante 
$$(\omega_{\C/S}.\omega_{\C/S}) \le \left ( 3 {g^2+\varepsilon \over 2g+1}-1\right ) \left ( \sum_{{\mathfrak p} \in S} \delta_{\mathfrak p} log N_{\mathfrak p}+\sum_{\sigma:K \rightarrow \CC} \delta_\sigma(\C_{\sigma})\right )+c'[K:\QQ]$$ 
o\`u $N_{\mathfrak p}$ d\'esigne le cardinal du corps r\'esiduel de $S$ en $\mathfrak p$ et o\`u $(\omega_{\C/S}.\omega_{\C/S})$ (autointersection du fibr\'e dualisant au sens d'Arakelov) et $\delta_\sigma(\C_{\sigma})$ ($\delta$ invariant de $\C$) sont d\'efinis \`a l'aide de la m\'etrique $\| \ \|_{Mod}$ sur $\omega_{\C/S}$.
\end{Cor}

\section{Courbes semi-stables hyperelliptiques}

\subsection{Rappels sur les torseurs de degr\'e $p$}

Soient $X$ un sch\'ema localement noetherien et $G$ un sch\'ema en groupes affine 
d'ordre $p$ sur $X$. Les $G$-torseurs au-dessus de $X$ sont alors classifi\'es \`a isomorphisme pr\`es par $\H^1(X, G)$ (cf \cite{Mil}, Corollary III.4.7 et Corollary III.2.10), la cohomologie \'etant prise au sens ppf (plate et de pr\'esentation finie). D'autre part, dans le cas o\`u $G$ est isomorphe \`a $\aalpha_p$ ou $\mmu_p$ et o\`u $X$ est une vari\'et\'e alg\'ebrique sur un corps parfait de caract\'eristique $p$, on a une description des $G$-torseurs au-dessus de $X$ en terme de formes diff\'erentielles.

\begin{The} 
\label{FormeDiff}
Soit $X$ comme ci-dessus. Notons $\cartier$ l'op\'erateur de Cartier sur $\H^0(X, \Omega^1_X)$, on a
$$\H^1(X, \aalpha_p) = \{ \omega \in \H^0(X, \Omega^1_X) \ | \ \der \omega = 0, \cartier \omega=0\}$$
et 
$$\H^1(X, \mmu_p) = \{ \omega \in \H^0(X, \Omega^1_X) \ | \ \der \omega = 0, \cartier\omega=\omega\}.$$
\end{The}

\noindent{\it Preuve :} Voir \cite{Mil}, proposition III.4.14. \Eproof
\medskip

D'autre part, lorsque $X$ est un sch\'ema affine, on a une desciption explicite des torseurs sous les groupes $p$-cycliques.

\begin{The}
\label{StructAp}
Soit $X$ un sch\'ema sur $\FF_p$ tel que $\H^1(X, \O_X)=0$ alors $$\H^1(X, \aalpha_p)=\H^0(X, \GGa)/\H^0(X, \GGa)^p.$$
 En particulier, si $X=\Spec(A)$ alors 
$$\H^1(X, \aalpha_p)=A/A^p$$ et les $\aalpha_p$ torseurs au-dessus de $X$ sont de la forme $\Spec(B)$ avec $B = A[y] / (y^p-u)$, o\`u $u$ est un \'el\'ement de $A$ unique \`a addition par une puissance $p$-i\`eme pr\`es. La forme diff\'erentielle associ\'ee est alors $\omega=du$.
\end{The}

\noindent{\it Preuve :}  Voir \cite{Mil}, \S III.4, {\it Infinitesimal coverings}, p. 128. \Eproof

\begin{The}
\label{StructMup}
Soit $X$ un sch\'ema sur $\FF_p$ tel que $\Pic(X)=0$. Alors $$\H^1(X, \mmu_p)=\H^0(X, \GGm)/\H^0(X, \GGm)^p.$$
 En particulier, si $X=\Spec(A)$ est un sch\'ema affine, alors 
$$\H^1(X, \mmu_p)=\GGm(A)/\GGm(A)^p,$$ et les $\mmu_p$-torseurs au-dessus de $X$ sont de la forme $\Spec(B)$ avec $B = A[y] / (y^p-u)$ o\`u $u$ est une unit\'e de $A$, unique \`a multiplication par la puissance $p$-i\`eme d'une unit\'e pr\`es. La forme diff\'erentielle associ\'ee est alors $\omega={du \over u}$.
\end{The}

\noindent{\it Preuve :}  Voir \cite{Mil}, \S III.4, {\it Kummer theory}, p. 125. \Eproof

\begin{Rem}
 Soient $X=\Spec(A)$ un sch\'ema affine sur $\FF_p$ et $Y$ un $\mmu_p$-torseur au-dessus de $X$. Une repr\'esentation de $Y$ sous la forme $\Spec(A[w]/(w^p-u))$ \'etant choisie, on peut oublier la structure de $\mmu_p$-torseur et munir $Y$ d'une structure (non canonique) de $\aalpha_p$-torseur au-dessus de $X$ en consid\'erant le morphisme $\aalpha_p \times Y \rightarrow Y$ correspondant \`a l'homomorphisme d\'efini par 
\begin{equation}
\begin{array}{cccl}
\label{mup_to_alphap}
A[w]/(w^p-u) & \rightarrow & A[x]/(x^p)  \otimes_A A[w]/(w^p-u) &\\
w & \mapsto & x \otimes 1 + 1 \otimes w &.\\
\end{array}
\end{equation}
On v\'erifie sans peine que ceci d\'efinit une structure de $\aalpha_p$-torseur. Cette structure d\'epend du choix de $u$.
\end{Rem}

\begin{Rem}
\label{StructZp}
 Dans le cas des $\ZZ/p\ZZ$-torseurs, on a un th\'eor\`eme analogue (cf. \cite{Mil}, \S III.4). Les $\ZZ/p\ZZ$-torseurs au-dessus d'un sch\'ema affine $\Spec(A)$ sur $\FF_p$ sont de la forme $\Spec(B)$ avec $B = A[y] / (y^p-y-u)$ o\`u $u$ est un \'el\'ement de $A$, unique \`a addition pr\`es par un \'el\'ement de la forme $a^p-a$ avec $a \in A$. \medskip
\end{Rem}

Dans la suite, nous aurons besoin d'une forme ``canonique'' pour les torseurs au dessus du disque formel, ceci nous est fournit par la proposition suivante.
 
\begin{Pro}
\label{FormeNormale}
Soient $k$ un corps alg\'ebriquement clos de caract\'eristique $p > 0$, $A=k((x))$, $B$ une extension de degr\'e $p$ de $A$ telle que $\Spec \ B$ soit un torseur sur $\Spec \ A$ sous $G$. Alors il existe un g\'en\'erateur $x'$ de $k[[x]]$, et un entier $m \in \ZZ$ tels que
\begin{enumerate}[(1)]
\item $B=A[y]/(y^p-y-{x'}^{-m})$ avec $m > 0$ et $p \nmid m$ si $G$ est \'etale. L'entier $m$ est alors d\'etermin\'e de fa\c con unique.
\item $B=A[y]/(y^p-{x'}^{-m})$ avec $p \nmid m$ si $G \cong \aalpha_p$. L'entier $m$ est alors d\'etermin\'e de fa\c con unique.
\item $B=A[y]/(y^p-{x'}^{-m})$ avec $p \nmid m$ ou $B=A[y]/(y^p-(1+x'Q))$ avec $Q \in k[[x']]$ si $G \cong \mmu_p$.
\end{enumerate}
\end{Pro}

\noindent{\it Preuve :} C'est une cons\'equence des th\'eor\`emes \ref{StructMup}, \ref{StructAp} et de la remarque \ref{StructZp}.  \Eproof

\begin{Def}
Avec les notations de la proposition pr\'ec\'edente, l'entier $m$ est appel\'e le {\it conducteur} de $B$ sur $A$ dans les cas $(1)$ et $(2)$.
\end{Def}

\if not
Ainsi, si $Y$ est un $\aalpha_p$-torseur au-dessus de $X$, il lui correspond une forme diff\'erentielle $\omega_Y$ sur $X$. Le fait que $\cartier \omega_Y$ soit nul impose que $\omega_Y$ est localement (pour la topologie ppf) de la forme $\der u$. De m\^eme, si $Y$ est un $\mmu_p$-torseur, on lui associe une forme diff\'erentielle $\omega_Y$ qui est localement (pour la topologie ppf) de la forme $\der u/u$ (car $\cartier\omega_Y = \omega_Y$). Par exemples, un $\aalpha_p$-torseur au dessus de $\AA^1_{\FF_p}$ sera d\'efini par une \'equation du type $y^p-f$ avec $f \in \FF_p[x]$. La forme diff\'erentielle associ\'ee sera alors $\der f$. De m\^eme, un $\mmu_p$-torseur au dessus de $\AA^1_{\FF_p}\setminus \{ 0 \}$ sera d\'efini par une \'equation du type $y^p-f$ avec $f \in \FF_p[x, 1/x]^*$ et la forme diff\'erentielle associ\'ee sera ${\der f \over f}$.
\fi

\subsection{D\'eformation de $p$-groupe cyclique en \'egale caract\'eri\-stique $p$}

On fixe un anneau de valuation discr\`ete $R$ d'\'egale caract\'eristique $p >0$.
 
Le cas de $\mmu_p$ n'est pas trait\'e ici car on ne peut pas d\'eformer $\mmu_p$ en $\ZZ/p\ZZ$ en \'egale caract\'eristique $p$ (ceci peut se voir en passant au dual de Cartier).

Soient $n \ge 0$ un entier et $t$ une uniformisante de $R$. Consid\'erons le morphisme de sch\'emas en groupes au-dessus de $\Spec(R)$ :
$$\phi_n : \GGa  \rightarrow \GGa$$ d\'eduit de l'homomorphisme
$$\begin{array}{llll}
\phi_n^\#: & R[x] & \rightarrow & R[x] \\
& x & \mapsto&  x^p-t^{(p-1)n} x .
\end{array}$$
Par convention, $t^0=1$. Notons $\M^n$ le noyau de $\phi_n$, on a 
\begin{equation}
\label{IsomAffine}
\M^n = \Spec \left (R[w]/(w^p-t^{(p-1)n}w) \right ).
\end{equation}
On a alors une suite exacte pour la topologie fppf
\begin{equation}\label{phi-n}
\xymatrix{0 \ar[r] &  \M^n \ar[r] & \GGa \ar[r]^{\phi_n} & \GGa  \ar[r] & 0}.
\end{equation} 
En effet, par d\'efinition de $\M^n$, il suffit de montrer que $\phi_n$ est surjective (pour la topologie fppf) et pour cela, il est suffisant de voir que $\M^n$ est plat sur $R$ (cf \cite{Mil}, exercise $2.19$), ce qui d\'ecoule de (\ref{IsomAffine}). 

\medskip
Si $n=0$, $\M^n$ est isomorphe \`a $\ZZ/p\ZZ$. Si $n>0$, la fibre g\'en\'erique de $\M^n$ est isomorphe \`a $\ZZ/p\ZZ$ et sa fibre sp\'eciale \`a $\aalpha_p$.

\medskip
\begin{Rem}
 Si $n>0$, la r\'eduction de la suite (\ref{phi-n}) modulo $t$
donne la suite sur $k$
$$\xymatrix{0 \ar[r] &  \aalpha_p \ar[r] & \GGa \ar[r]^{x \mapsto x^p} & \GGa  \ar[r] & 0}$$
dont l'exactitude est connue. Dans le cas $n=0$, on retrouve la suite exacte d'Artin-Schreier.
\end{Rem}

\begin{Pro}
\label{Relevement_torseur}
Soient $X$ un sch\'ema de type fini sur $\Spec(R)$ et $n\ge 0$. On suppose que $\H^1_{Zar}(X_s, \O_{X_s})=0$ et que le morphisme $\H^0_{Zar}(X, \O_X) \rightarrow \H^0_{Zar}(X_s, \O_{X_s})$ est surjectif.
Alors tout $\M^n_s$ torseur sur $X_s$ se rel\`eve en un $\M^n$-torseur sur $X$. 
\end{Pro}

\noindent{\it Preuve :}  Il suffit de montrer que le morphisme de sp\'ecialisation $$\H^1(X, \M^n) \rightarrow \H^1(X_s, \M^n_s))$$ est surjectif. Pour cela, consid\'erons la suite exacte \eqref{phi-n}.
On a un diagramme commutatif \`a lignes exactes donn\'e par la suite exacte longue de cohomologie (fppf) :
$$\xymatrix{
\H^0(X, \GGa) \ar[r]\ar[d] & \H^1(X, \M^n) \ar[r]\ar[d] & \H^1(X, \GGa) \ar[d] \\
\H^0(X_s, \GGa) \ar[r] & \H^1(X_s, \M^n_s) \ar[r] & \H^1(X_s, \GGa).}$$
Pour tout $i \ge 0$ et tout sch\'ema $Y$ on a (\cite{Mil}, proposition III.3.7) 
$$\H^i_{fppf}(Y, \GGa)=\H^i_{fppf}(Y, \O_Y)=\H^i_{Zar}(Y, \O_Y).$$
Les hypoth\`eses sur $X$ montrent qu'on a un diagramme commutatif
$$\xymatrix{
\H^0(X, \GGa) \ar[r]\ar@{->>}[d] & \H^1(X, \M^n) \ar[d]  \\
\H^0(X_s, \GGa) \ar@{->>}[r] & \H^1(X_s, \M^n_s).}$$
Ce qui prouve que le morphisme $\H^1(X, \M^n) \rightarrow \H^1(X_s, \M^n_s)$ est surjectif. \Eproof

\begin{Cor}
\label{ReleveTorseur}
Si $X$ est affine, alors tout $\M^n_s$-torseur sur $X_s$ se rel\`eve en un $\M^n$-torseur sur $X$. 
\end{Cor}

\begin{Rem}
\label{StructMn}
Si $X=\Spec(A)$ est affine, la suite exacte (\ref{phi-n}) montre que $$\H^1(X, \M^n)=A/\phi_n(A)$$
avec $\phi_n:a \mapsto a^p-t^{(p-1)n}a$.
 Par suite, les $\M^n$-torseurs au dessus de $X$ sont affines de la forme $$\Spec(A[w]/(w^p-t^{(p-1)n} w-a))$$ o\`u $a$ est un \'el\'ement de $A$ unique \`a l'addition d'un \'el\'ement de $\phi_n(A)$ pr\`es.
En particulier, \'etant donn\'e le th\'eor\`eme \ref{StructAp} et la remarque \ref{StructZp}, on voit que pour relever des $\M^n$-torseurs au-dessus d'un sch\'ema affine, il suffit de relever une \'equation.
\end{Rem}

\subsection{Rev\^etement $p$-admissible et courbes hyperelliptiques}

\begin{Def}
Nous appelons courbe {\it courbe semi-stable} sur un corps alg\'ebrique\-ment clos $k$ 
une courbe connexe, propre, r\'eduite de dimension 1 sur $\Spec(k)$ 
dont les composantes irr\'eductibles sont lisses et dont les points singuliers sont doubles ordinaires.
\end{Def}

Soit $C$ une courbe semi-stable sur un corps alg\'ebriquement clos, on lui associe un graphe appel\'e {\it graphe d'intersection de $C$} dont les sommets sont les composantes irr\'eductibles de $C$ et les ar\^etes sont les points doubles, deux sommets \'etant reli\'es par les ar\^etes correspondant aux points doubles de l'intersection des deux composantes.
La plupart des graphes consid\'er\'es seront orient\'es, cette orientation \'etant fix\'ee arbitrairement.

\medskip
Pour toute vari\'et\'e alg\'ebrique $V$ sur $k$, on note $V_{lisse}$ l'ouvert de lissit\'e de $V$. Soit $Y$ une courbe semi-stable sur $k$. On note $\Gamma_Y$ le graphe d'intersection de $Y$. Pour tout sommet $v$ de $\Gamma_Y$, on note $F_v$ la composante irréductible de $Y$ correspondante à $v$ et $F_v^0:=F_v\cap Y_{lisse}$. 

\begin{Def} 
\label{DonneeCompatible}
Soit $k$ un corps alg\'ebriquement clos de caract\'eristique $p > 0$. Soit $f : Y\to X$ un morphisme fini de courbes projectives lisses sur $k$. 
Soit $U$ un ouvert affine de $X$ tel que $f^{-1}(U)\to U$ soit un torseur
sous un des groupes alg\'ebriques $G=\ZZ/p\ZZ$, $\aalpha_p$ et 
$\mmu_p$ sur $k$. La structure du torseur est induit par un \'el\'ement $u\in\O_X(U)$
\`a la mani\`ere des th\'eor\`emes \ref{StructAp}, \ref{StructMup} et de la remarque \ref{StructZp}. 
Soit $P\in Y$ un point ramifi\'e de $f$. Nous allons lui 
associer des entiers $m(P)$ ({\it conducteur}) et $h(P)$ ({\it r\'esidu}). 
Notons $Q=f(P)$. Si $G=\ZZ/p\ZZ$, on pose $m(P)$ \'egal au 
conducteur de Hasse de l'extension $\O_{X, Q}\to \O_{Y, P}$, et $h(P)=0$. 
Si $G=\aalpha_p$, on pose 
$$m(P)=-(1+\Ord_Q du), \quad h(P)=0.$$ 
Enfin, si $G=\mmu_p$ (donc $u\in\O_X(U)^*$), on pose 
$$m(P)=-(1+\Ord_Q du/u), \quad h(P)=\Res_Q du/u.$$ 
Ces entiers ne d\'ependent pas du choix de $u$, mais uniquement de la 
structure de torseur au-dessus de $U$.
\end{Def}

\begin{Lem}
\label{NN} Conservons les notations ci-dessus. 
\begin{enumerate}[\rm (a)]
\item Si $G=\aalpha_p$ ou $\ZZ/p\ZZ$, alors $m(P)\ne 0$.
\item Si $G=\mmu_p$, alors $m(P)=0$ si et seulement si $n:=\Ord_P u$ est 
premier \`a $p$, et on a alors $h(P)=n$ (dans $k$). 
\item Supposons $G=\mmu_p$ et $X=\PP^1_k$. Alors quitte \`a multiplier $u$ 
par une puissance $p$-i\`eme dans $\O_X(U)^*$ (ce qui ne change pas la 
structure de torseur), on a $u\in \O_{X,f(P)}^*$ pour tout 
$P\in Y$ tel que $m(P)\ne 0$. \end{enumerate}
\end{Lem}

\noindent{\it Preuve}. (a)--(b) Si $G$ est \'etale, on a $m(P)>0$ par la 
d\'efinition du conducteur de Hasse. Dans le cas $G$ radiciel, les propri\'et\'es 
se v\'erifient ais\'ement en \'ecrivant $u$ sous la forme $u=\pi^n u_0$ avec $\pi$ 
un param\`etre local de $X$ en $Q$ et $u_0\in\O_{X, Q}^*$. 

(c) On a $(u)=pD_1+D_0$ o\`u les $D_i$ sont des diviseurs sur $X$, et 
$D_0$ est \`a coefficients premiers \`a $p$. Soit $Q_0$ un point du support
de $D_0$ (qui est non vide car $u$ serait une puissance $p$-i\`eme sinon).
Alors il existe $w\in k(X)$ tel que $(w)=D_1-(\deg D_1)Q_0$. Il suit que
$w\in \O_X(U)^*$ et que $(u/w^p)$ est un diviseur dont les coefficients
(non nuls) sont premiers \`a $p$. Si $m(P)\ne 0$, alors $\Ord_P(u/w^p)$ est 
divisible par $p$, donc nul. \Eproof

\begin{Def}
\label{RevpAdmissible}
Soit $k$ un corps alg\'ebriquement clos de caract\'eristique $p>0$.
Un {\it rev\^etement $p$-admissible} de courbes 
semi-stables consiste en la donn\'ee de: 
\begin{enumerate}[\rm (1)]

\item un morphisme fini surjectif $f : Y\to X$ de courbes semi-stables
sur $k$;

\item une application qui \`a chaque sommet $v$ de $\Gamma_Y$ associe un
couple $(G_v, u_v)$, o\`u $G_v$ est l'un des groupes alg\'ebriques 
$\ZZ/p\ZZ$, $\aalpha_p$ et $\mmu_p$ sur $k$, et o\`u $u_v$
est une fonction rationelle sur $f(F_v^0)$. 
\end{enumerate}
Ces donn\'ees doivent v\'erifier les propri\'et\'es suivantes:
\begin{enumerate}[\rm (a)]
\item $f^{-1}(X_{lisse})=Y_{lisse}$;

\item pour tout sommet $v$ de $\Gamma_Y$, il existe un ouvert non vide $U_v \subset f(F_v^0)$ tel que $f^{-1}(U_v)\to U_v$ 
soit un torseur sous $G_v$, d\'efini par $u_v$ \`a la mani\`ere des th\'eor\`emes 
\ref{StructAp}, \ref{StructMup} et de la remarque \ref{StructZp} ($u_v=0$ si $f^{-1}(U_v)$ n'est pas
connexe). En particulier, $f^{-1}(f(F_v))\to f(F_v)$ est un morphisme fini de
degr\'e $p$. On dira que $v$ est {\it \'etale} (resp. {\it radiciel})
si $G_v$ est \'etale (resp. radiciel). 

\item Soit $a$ une ar\^ete (orient\'ee) de $\Gamma_Y$ correspondant \`a un point
d'intersection $P$. Supposons que $f^{-1}(f(P))=\{ P \}$. On dira que 
$a$ est une {\it ar\^ete ramifi\'ee}. 
Alors on demande que l'une des extr\^emit\'es 
de $a$ soit radicielle. Soit $v$ le sommet d'origine de $a$. Posons 
$m(a)=m(P)$ et $h(a)=h(P)$. Alors on impose la {\it condition de 
compatibilit\'e}:
$$m(a)=-m(\bar{a}), \quad h(a)=-h(\bar{a}), $$
o\`u $\bar{a}$ est l'ar\^ete oppos\'ee de $a$. 
\end{enumerate}
On notera un rev\^etement $p$-admissible par $(f : Y\to X, \{ (G_v, u_v) \}_v)$.
On dira que $(f : Y\to X, \{ (G_v, u'_v) \}_v)$ est {\it \'equivalent} \`a 
$(f : Y\to X, \{ (G_v, u_v) \}_v)$ si $u_v$ et $u'_v$ induisent la m\^eme 
structure de torseur au-dessus de $U_v \subset f(F_v^0)$. 
\end{Def}

Dans toute la suite, on prendra pour $U_v$ le plus grand ouvert de $f(F_v^0)$ tel que $f^{-1}(U_v) \rightarrow U_v$ soit un torseur sous $G_v$ d\'efini par $u_v$.

\begin{Def}[{\rm \cite{HM}, \S 4}]
Soient $k'$ un corps alg\'ebriquement clos de caract\'eristique {\bf diff\'erente} de $p$, 
$Y$ et $X$ deux courbes semi-stables sur $k'$ et $f:Y \rightarrow X$ un morphisme fini. On dit que $f$ est un {\it rev\^etement $p$-admissible} si les conditions suivantes sont r\'ealis\'ees :
\begin{enumerate}[-]
\item il existe une action de $\ZZ/p\ZZ$ sur $Y$ libre sur un ouvert dense telle que $f$ s'identifie au morphisme quotient $Y \rightarrow Y/(\ZZ/p\ZZ)$.
\item soit $P$ un point fixe de $Y$ sous l'action de $\ZZ/p\ZZ$, alors l'action de $\ZZ/p\ZZ$ sur $\O_{Y, P}$ est kummerienne (c'est \`a dire que si on note $\chi_1$ et $\chi_2$ les caract\`eres $\ZZ/p\ZZ \rightarrow k^*$ d\'eduits de l'action de $\ZZ/p\ZZ$ sur chaque branche de $Y$ en $P$, on a $\chi_1 \chi_2=1$). 
\end{enumerate} 
\end{Def}

\begin{Def}
\label{CourbeStableHyperelliptique}
Soient $k$ un corps et $C \rightarrow \Spec(k)$ une courbe stable de genre arithm\'etique $g \ge 2$. On dira que $C$ est {\it hyperelliptique} s'il existe un anneau de valuation discr\`ete $R_1$ dont le corps r\'esiduel est une extension finie  $k_1$ de $k$, une courbe $\C \rightarrow \Spec \ R_1$ stable dont la fibre g\'en\'erique est lisse et hyperelliptique, et dont la fibre sp\'eciale est isomorphe \`a $C_{k_1}$.
\end{Def}

\begin{Rem}
 Soient $k$ un corps, $C \rightarrow \Spec(k)$ une courbe lisse hyperelliptique de genre $g \ge 2$, alors $C$ est hyperelliptique au sens habituel. 
\end{Rem}

\begin{Pro} 
\label{Reduction_Admissible}
Supposons que $K$ est de caract\'eristique z\'ero. Soit $\C\rightarrow \Spec(R)$ une courbe stable dont la fibre g\'en\'erique est lisse et munie d'une action non triviale d'un groupe $G$ $p$-cyclique. Alors apr\`es extension finie de $R$, il existe un \'eclatement $\tilde \C \rightarrow \C$ et, si $\car \ k=p$, des donn\'ees $\{(G_v, u_v)\}_{v\in\Gamma_{\C_s}}$ tels que $\tilde \C_s \rightarrow (\tilde \C/G)_s$ (resp. $(\tilde \C_s \rightarrow (\tilde \C/G)_s, \{(G_v, u_v)\}_{v\in\Gamma_{\C_s}})$) soit un rev\^etement $p$-admissible d'une courbe semi-stable si $\car \ k\not = p$ (resp. $\car \ k = p$).
\label{Descente}
\end{Pro}

\noindent{\it Preuve :} On peut supposer $k$ alg\'ebriquement clos.
Supposons dans un premier temps que $k$ est de caract\'eristique diff\'erente de $p$. Quitte \`a faire une extension de $R$, on peut supposer que les \'epaisseurs des points singuliers sont paires. Notons $\tilde \C$ la courbe obtenue de $\C$ en \'eclatant les points $P$ de la fibre sp\'eciale pour lesquels $G$ ne laisse pas fixe les branches analytique en $P$. Par suite, \cite{Sai} Th\'eor\`eme 5.7 montre que $\tilde \C_s \rightarrow (\tilde \C/G)_s$ est $p$-admissible.

Supposons d\'esormais que $k$ est de caract\'eristique $p$. Quitte \`a augmenter $K$, on peut supposer que tout les points fixes de $\C_\eta$ sous $G$ sont rationnels. Consid\'erons le mod\`ele minimal $\tilde\C$ de $\C$ d\'eployant les sp\'ecialisations des points de la fibre g\'en\'erique fixes sous $G$ en des points lisses et distincts. Par unicit\'e du mod\`ele minimal $\tilde \C$, l'action de $G$ se prolonge \`a $\tilde \C$. Quitte \`a faire une extension ramifi\'ee de l'anneau $R$, on peut supposer les \'epaisseurs des points doubles paires. Soit $P \in \tilde \C_s$ un point singulier fixe sous $G$  et tel que $G$ permutte les branches analytiques en $P$. Soit $f_1:\tilde \C_1 \rightarrow \tilde \C$ l'\'eclatement de $\tilde \C$ en $P$. Alors $G$ op\`ere sur $\tilde \C_1$ et les points singuliers de $\tilde \C_{1, s}$ contenus dans $f_1^{-1}(P)$ ne sont pas fixes par $G$. On peut donc supposer que $G$ ne permutte les branches analytiques en aucun point de $\tilde \C$. Soient $C_1$ et $C_2$ deux composantes \'etales de $\tilde \C_s$ et $P \in C_1 \cap C_2$. Soit $f_2:\tilde \C_2 \rightarrow \tilde \C$ l'\'eclatement de $\tilde \C$ en $P$. Alors $G$ agit sur $\tilde \C_2$ et l'intersection des transform\'ees strictes de $C_1$ et $C_2$ dans $\tilde \C_2$ ne contient aucun point de $f_2^{-1}(P)$. On peut donc supposer qu'aucun point fixe de $\tilde \C_s$ n'est intersection de deux composantes \'etales.

Pour chaque sommet $v$ de $\Gamma_{\C_s}$, la r\'eduction de l'action de $G$ induit une structure de $G_v$-torseur ($G_v$ est isomorphe \`a $\aalpha_p$, $\mmu_p$ ou $\ZZ/p\ZZ$) sur un ouvert dense de la composante correspondant \`a $v$ (d'apr\`es l'\'etude faite dans dans \cite{Hen_p_adique}, Proposition 1.6 sur la r\'eduction des $\mmu_p$-torseur) et on peut donc lui associer $u_v$ \`a la mani\`ere des th\'eor\`emes \ref{StructAp}, \ref{StructMup} et de la remarque \ref{StructZp}. Nous allons montrer que les donn\'ees $(\tilde \C_s \rightarrow (\tilde\C/G)_s, \{(G_v, u_v)\})$ d\'efinissent un rev\^etement $p$-admissible, c'est \`a dire que les conditions $(a)$, $(b)$ et $(c)$ de la d\'efinition \ref{RevpAdmissible} sont satisfaites :

\begin{enumerate}[(a)]
\item Comme $G$ n'\'echange pas les branches d'un point singulier, l'image r\'eciproque du lieu lisse de $(\tilde \C/G)_s$ est le lieu lisse de $\tilde \C_s$.

\item Le fait que $f$ soit un $G_v$-torseurs en restriction \`a un ouvert non vide de chaque composante (correspondant \`a $v$) du lieu lisse provient de la construction.

\item Les compatibilit\'es au niveau des conducteurs et des r\'esidus proviennent de \cite{Hen_p_adique}, Proposition 3.1. \Eproof
\end{enumerate}

\section{Un th\'eor\`eme de rel\`evement en \'egale caract\'e\-ristique}
\label{thm-rel}

On d\'emontre dans cette section un th\'eor\`eme de rel\`evement des rev\^etements $p$-admissibles en \'egale caract\'eristique $p$ (Th\'eor\`eme \ref{Relevement}). Pour cela, on commence par donner quelques r\'esultats sur le rel\`evement des $\ZZ/p\ZZ$ et $\aalpha_p$-torseurs en \'egale caract\'eristique.

\subsection{Graphe de Hurwitz}\label{arb-h}

Cette notion est tir\'ee de \cite{Hen_p_adique}, d\'efinition 2.1. Nous \'etendons toutefois la notion au cas o\`u $\car K=p$. On prendra pour convention $\nu(p)= \infty$ si $K$ est un corps de caract\'eristique $p$.

Soit $\Gamma$ un graphe orient\'e, nous noterons ${\bf Som} \Gamma$ (resp. ${\bf Ar} \Gamma$, resp. ${\bf Ar}^+ \Gamma$ ) l'ensemble des sommets (resp. des ar\^etes orient\'ees, resp. des ar\^etes positives) de $\Gamma$. Si $G$ est un groupe abstrait fini agissant sur $\Gamma$, on d\'esignera par ${\bf Som}^G \Gamma$ (resp. ${\bf Ar} ^G\Gamma$) l'ensemble des sommets (resp. des ar\^etes orient\'ees) de $\Gamma$ fixes sous l'action de $G$. Si $v$ est un sommet de $\Gamma$ fixe sous $G$, ${\bf Ar}^+ (s)$ sera l'ensemble des ar\^etes orient\'ees fixe sous $G$ et ayant pour origine $s$.

Soit $G$ un groupe $p$-cyclique agissant sans inversion (c'est \`a dire que pour toute ar\^ete $a$ et pour tout $\tau \in G$ on a $\tau a \not = \bar a$) sur $\Gamma$. Consid\'erons la donn\'ee ${\cal H}=(\varepsilon, g, d, m, h)$, o\`u :
\begin{enumerate}[-]
\item $\varepsilon : {\bf Ar}\Gamma \rightarrow \NN$ est une application dans l'ensemble des entiers $\ge 0$. Si $v$ est un sommet, $\varepsilon(v)$ sera appel\'e l'{\it \'epaisseur} de $s$.
\item $g:{\bf Som}\Gamma \rightarrow \NN$ est une application. Si $v$ est un sommet, on appellera $g(v)$ le {\it genre} de $v$.
\item $d:{\bf Som}\Gamma\rightarrow (p-1)\NN \cap [0, \nu(p)] $ est une application; pour un sommet $v$ fix\'e, l'entier $d(v)$ sera appel\'e sa {\it diff\'erente}.
\item $m:{\bf Ar}^G\Gamma \rightarrow \ZZ$ est une application; pour une ar\^ete orient\'ee $a$ fixe par $G$, l'entier $m(a)$ sera appel\'e {\it conducteur} de $a$.
\item $h:{\bf Ar}^G\Gamma \rightarrow \ZZ/p\ZZ$ est une application; pour une ar\^ete orient\'ee $a$ fixe par $G$, l'entier $h(a)$ sera appel\'e {\it r\'esidu} de l'ar\^ete.
\end{enumerate}

Un sommet $v$ est dit {\it additif} si $\nu(p) > d(v) > 0$, {\it \'etale} si $d(v)=0$ et, si $K$ est de caract\'eristique diff\'erente de $p$, on dira que $v$ est {\it multiplicatif} si $d(v)=\nu(p)$.        
La donn\'ee $\cal H$ est appel\'ee $K${\it -donn\'ee de Hurwitz} sur $\Gamma$, relativement \`a l'action de $G$, si de plus les conditions suivantes sont satisfaites :
\begin{enumerate}[-]
\item L'action de $G$ est compatible avec l'\'epaisseur, le genre et la diff\'erente, c'est-\`a-dire que pour tout $\tau$ dans $G$, pour tout sommet $v$ et toute ar\^ete $a$ on a $\varepsilon(\tau.a)=\varepsilon(a)$, $g(\tau.v)=g(v)$
 et $d(\tau.v)=d(v)$.
\item Pour toute ar\^ete $a$ fixe par $G$, si $m(a)$ est non nul alors $m(a)$ est premier \`a $p$.

\item Pour tout sommet $v$ \'etale fixe par $G$ et pour toute ar\^ete orient\'ee $a$ d'origine $v$ fixe par $G$, on a $m(a) > 0$ et il existe un entier $g'(v)\ge 0$ tel que

\noindent $2g(v)-2=p(2g'(v)-2)+\sum_{a \in {\bf Ar}^G(v)} (m(a)+1)(p-1)$ \hfill (``formule de Hurwitz'')
\item Pour tout sommet $v$ radiciel, on a les relations

\noindent $\sum_{a \in {\bf Ar}^G(v)} (m(a)+1)=2-2g(v)$ \hfill (``degr\'e du diviseur canonique'')

\noindent $\sum_{a \in {\bf Ar}^G(v)} h(a)=0$ \hfill (``th\'eor\`eme des r\'esidus'')   
\item Pour toute ar\^ete $a$ fixe par $G$, on a $m(\bar a)=-m(a)$, $h(\bar a)=-h(a)$ et 

\noindent $d(t(a))-d(o(a))=m(a) \varepsilon(a)(p-1)$ \hfill (``loi de variation de la diff\'erente'')
\end{enumerate}

\medskip

\begin{Ex} 
\label{ExempleReduction}
Soit $R$ un anneau de valuation discr\`ete complet  d'in\'egale caract\'eri\-stique $(0, p)$. La donn\'ee d'une courbe $\pi:\C \rightarrow \Spec(R)$ stable dont la fibre g\'en\'erique est lisse, d'une action d'un groupe $p$-cyclique $G$ sur $\C$ fournit naturellement un graphe de Hurwitz. En effet, notons $\tilde \C$ le mod\`ele minimal d\'eployant en des points lisses distincts les sp\'ecialisations des points de la fibre g\'en\'erique de $\C$ qui sont fixe sous $G$. Notons $\Gamma$ le graphe dont l'ensemble des sommets est compos\'e des composantes irr\'eductibles de $\tilde \C_s$ et des points fixes de $\C_\eta$. Deux sommets $v$ et $w$ sont reli\'es par une ar\^etes lorsque
\begin{enumerate}[-]
\item $v$ et $w$ correspondent \`a deux composantes irr\'eductibles de $\tilde \C_s$ qui se coupe en un point, il y a alors autant d'ar\^etes joignant $v$ et $w$ que de points d'intersection entre ces deux composantes.
\item $v$ correspond \`a une composante de $\tilde \C_s$ et $w$ \`a un point de $\tilde \C_\eta$ se sp\'ecialisant sur cette composante. 
\end{enumerate}
Pour plus de d\'etail, voir \cite{Hen_p_adique}, Proposition 3.1. Nous nous contentons ici d'illustrer ce fait g\'en\'eral par un exemple. Supposons que $\C$ soit une courbe hyperelliptique telle que $\tilde \C_s$ soit donn\'ee par la figure \ref{exemple}.

\begin{figure}[h]
\centerline{\input{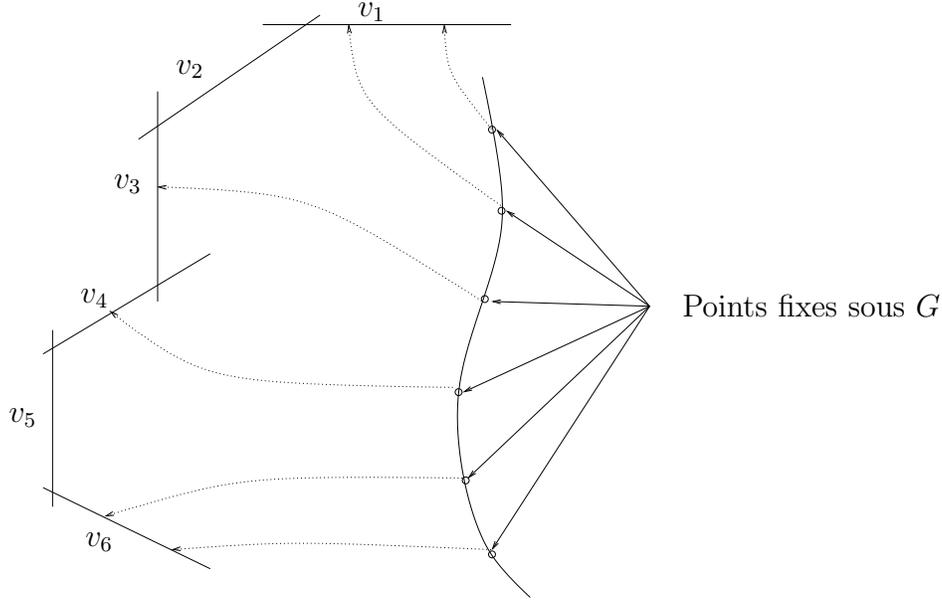}\hskip 3cm}
\caption{Fibre sp\'eciale de $\tilde \C$ et sp\'ecialisation des points fixes}
\label{exemple}
\end{figure}

Notons $\Gamma$ le graphe obtenu comme pr\'ec\'edement.
Pour chaque sommet $v_i$ du graphe de $\Gamma$ correspondant \`a une composante irr\'eductible de $\tilde \C_s$, la proposition 1.6 de \cite{Hen_p_adique}, il existe un ouvert non vide de $F_{v_i}^0$ naturellement muni d'une structure de $\aalpha_2$, $\mmu_2$ o\`u $\ZZ/2\ZZ$-torseur (en fonction de la valuation de la diff\'erente). De plus, on peut associer \`a $v_i$ le genre de la composante correspondante et la diff\'erente de l'anneau local de son point g\'en\'erique (dans $\C$). Pour chaque sommet $v_i$ correspondant \`a un point fixe sous $G$ se sp\'ecialisant dans une composant $v$, le genre et la diff\'erente sont fix\'es \`a z\'ero. Pour chaque ar\^ete $a$ de $\Gamma_{\C_s}$, on peut donc associer un conducteur et un r\'esidu (comme dans la d\'efinition \ref{NN}), une \'epaisseur (l'\'epaisseur du point double correspondant dans $\C$).

L'exemple suivant est alors un arbre de Hurwitz (les sommets provenant de composantes irr\'eductibles sont repr\'esent\'es par des points et les sommets provenant de point fixe par des croix).

\centerline{\input{hurwitz.pstex_t}}

\end{Ex}

\subsection{\'Enonc\'e du th\'eor\`eme de r\'ealisation}

Dans toute cette section, on supposera que les sch\'emas en groupes intervenant dans la d\'efinition des rev\^etements $p$-admissibles sont isomorphes soit \`a $\ZZ/p\ZZ$, soit \`a $\aalpha_p$. Cette hypoth\`ese suppl\'ementaire sera justifi\'ee lors de la d\'emonstration du th\'eor\`eme \ref{ChangeTorseur}. On suppose de plus dans cette section que $\car \ K=p$.

\begin{Def} \label{DefAdapt}
Soient $R$ un anneau de valuation discr\`ete complet d'\'egale caract\'eristique $p$, $\X$ une courbe formelle semi-stable sur $R$ (c'est \`a dire que la fibre sp\'eciale de $\X$ est semi-stable), de fibre sp\'eciale $\X_s$ et $(f_s:\Y_s \rightarrow \X_s, \{(G_v, u_v)\}_v)$ un rev\^etement $p$-admissible. Un graphe de Hurwitz $\Gamma$ sera dit {\it adapt\'e} \`a $\X$ et au rev\^etement $p$-amissible $(f_s, \{G_v, u_v\}_v)$ si 
\begin{enumerate}[-]
\item l'ensemble des sommets de $\Gamma$ est compos\'e des sommets de $\Gamma_{\Y_s}$ et des points r\'eguliers de $\Y_s$ en lesquels $f_s$ n'est pas un torseur (i.e. les points de $F_{v}^0 \setminus f^{-1}(U_v)$).
\item le genre des sommets de $\Gamma$ correspondant \`a une composante irr\'eductible de $\Y_s$ est \'egal au genre de cette composante, et nul sinon.
\item l'ensemble des ar\^etes de $\Gamma$ est compos\'e des ar\^etes de $\Gamma_{\Y_s}$ et des couples $(v, P)$ o\`u $v$ est une composantes irr\'eductibles de $\Y_s$, $P$ est un point de $F_v^0$  en lequel $f_s$ n'est pas un torseur.
\item les conducteurs et r\'esidus en les ar\^etes de $\Gamma$ sont \'egaux aux conducteurs et r\'esidus en les points de $\Y_s$ (comme dans la d\'efinition \ref{DonneeCompatible}) correspondant.
\item pour toute ar\^ete $a$ de $\Gamma$ correspondant \`a un point singulier $P$ de $\Y_s$, $p\varepsilon(a)$ est \'egal \`a l'\'epaisseur de $f_s(P)$ dans $\X$ (en particulier, cela impose que les \'epaisseurs des points singuliers de $\X_s$ soient multiples de $p$).
\end{enumerate}
\end{Def}

\begin{The}
Soient $G$ un groupe $p$-cyclique abstrait, $R$ un anneau de valuation discr\`ete complet, $\X$ une courbe formelle semi-stable sur $R$, de fibre sp\'eciale $\X_s$ et $(f_s:\Y_s \rightarrow \X_s, \{(G_v, u_v)\}_v)$ un rev\^etement $p$-admissible tel que $f_s$ ne soit pas un hom\'eomorphisme. On fait les hypoth\`eses suivantes :
\begin{enumerate}[{\rm (\rm H}$_1${\rm )}]
\item pour tout $v$, $G_v=\ZZ/p\ZZ$ ou $\aalpha_p$. 
\item soient $v$ un sommet de $\Gamma_{\Y_s}$ et $P$ un point de $F_v^0$ en lequel $f_s$ n'est pas un torseur, alors $G_v \cong \aalpha_p$ et l'ordre de $u_v$ en $P$ est de la forme $1-pn$ avec $n \ge 0$.
\end{enumerate}
Si $\Gamma$ est un graphe de Hurwitz adapt\'e \`a $(\X, (f_s, \{(G_v, u_v)\}_v))$ alors il existe une courbe formelle semi-stable $\Y$ sur $R$ munie d'une action de $G$ et un $R$ isomorphisme $\Y/G \rightarrow \X$ tel que la fibre sp\'eciale du morphisme $\Y \rightarrow \X$ s'identifie \`a $f_s$.
\label{Relevement}
\end{The}

\begin{Rem}
Ce th\'eor\`eme est un analogue en \'egale caract\'eristique de certains r\'esultats de \cite{Hen_p_adique} et de \cite{Sai1}, \cite{Sai2}, \cite{Sai3} en in\'egale caract\'eristique.
\end{Rem}

\begin{Cor} 
\label{RelevementCourbe}
Soient $Y$ et $X$ deux courbes semi-stables sur $\Spec(k)$ et $(f:Y \rightarrow X, \{(G_v, u_v)\}_v)$ un rev\^etement $p$-admissible comme dans le th\'eor\`eme \ref{Relevement}. Soient $\cal X$ un rel\`evement de $X$ sur $\Spec \ R$ dont les \'epaisseurs en chaque point singulier est un multiple de $p$ et $\Gamma$ un graphe de Hurwitz adapt\'e \`a ${\cal X}$ et \`a  $(f, \{(G_v, u_v)\})$. Alors il existe une courbe semi-stable $\cal Y$ sur $\Spec(R)$ et un rev\^etement $p$-cyclique $f':\cal Y \rightarrow {\mathcal X}$ tels que la fibre g\'en\'erique de ce morphisme soit un rev\^etement galoisien de degr\'e $p$ entre deux courbes lisses, et que $f_s'=f$.
\end{Cor}

\noindent{\it Preuve :}  Il suffit d'alg\'ebriser les donn\'ees du th\'eor\`eme, ce qui est possible par le th\'eor\`eme d'ag\'ebrisation \cite{EGA}, Corollaire III.5.1.6. \Eproof
\medskip

	Pour la d\'emonstration de ce th\'eor\`eme, on aura besoin d'un th\'eor\`eme de recollement formel que nous rappelons ci-dessous.

	Soit $\X$ une courbe formelle semi-stable sur $R$, de fibre g\'en\'erique lisse et g\'eom\'etri\-quement connexe. Soient $x_1, \ldots, x_n$ des points ferm\'es lisses de la fibre sp\'eciale de $\X$ avec $x_i \in F_{v_i}^0$. On note $\Gamma$ le graphe orient\'e d'intersection de $\X$. Pour tout sommet $v$ de $\Gamma$, on notera $\X_v$ un ouvert formel de $\X$ dont la fibre sp\'eciale s'identifie \`a $F_v^0 \setminus\{x_1, \ldots, x_n\}$. Pour tout $i \in \{1, \ldots, n\}$, on notera $\X_i=\Spf \hat\O_{\X, x_i}$ et si $\eta_i$ est le point g\'en\'erique de $F_{v_i}$, on notera $\alpha_i:\Spf \widehat{({\hat\O_{\X, x_i}})_{\eta_i}} \rightarrow \X_i$ et $\beta_i:\Spf \widehat{({\hat\O_{\X, x_i}})_{\eta_i}} \rightarrow \X_{v_i}$. Pour toute ar\^ete $a$ de $\Gamma$ correspodant \`a un point double $P$, on notera $\X_a$ la couronne formelle $\Spf \ \hat \O_{\X, P}$ et $\X^o_a=\Spf \widehat {(\hat \O_{\X, P})_{\eta_{a, o}}}$, o\`u $\eta_{a, o}$ est le point g\'en\'erique de $\X_a$ correspondant au bord origine. Pour toute ar\^ete $a$ de $\Gamma$, on dispose alors des morphismes canoniques $\alpha_a:\X_a^o \rightarrow \X_a$ et $\beta_a:\X_a^o \rightarrow \X_{o(a)}$. On note $\Delta \X$ le sch\'ema formel $$\left(\coprod_{i \in \{1, \ldots, n\}} \X_i \right) \coprod\left (\coprod _{v \in {\bf Som}(\Gamma)} \X_v \right ) \coprod  \left (\coprod_{a \in {\bf Ar}(\Gamma)} \X_a \right )$$ et $q$ la projection naturelle $q:\Delta \X \rightarrow \X$.

	Soit $\F$ un module coh\'erent sur $\X$. Si $v$ est un sommet  (resp. $a$ est une ar\^ete, resp. $i \in \{1, \ldots, n\}$) de $\Gamma$, on note $\F_v$ (resp. $\F_a$, resp. $\F_i$) la restriction de $\F$ \`a $\X_v$ (resp. $\X_a$, resp. $\X_i$). Autrement dit $q^*\F=((\F_v), (\F_a), (\F_i))$. Pour tout $i\in\{1, \ldots, n\}$, on a alors des isomorphismes $\phi_i(\F):\alpha_{i}^* \F_i \rightarrow \beta_{i}^* \F_{v_i}$. Pour toute ar\^ete $a$ de $\Gamma$, on a alors un isomorphisme naturel $\phi_a(\F) : \alpha_{a}^* \F_a \rightarrow \beta_{a}^* \F_{o(a)}$. Plus g\'en\'eralement, si $\Delta \F = ((\F_v), (\F_a), (\F_i))$ est un module coh\'erent sur $\Delta \X$, on appelle donn\'ee de descente sur $\Delta \F$ une collection d'isomorphismes $\phi_a : \alpha_{a}^*\F_a \rightarrow \beta_{a}^* \F_{o(a)}$ et $\phi_i: \alpha_{i}^*\F_i \rightarrow \beta_{i}^* \F_{v_i}$. Les modules coh\'erents munis d'une donn\'ee de descente forment une cat\'egorie not\'ee $(\Delta {\mathcal Coh}_\X)$.

\begin{Pro}
\label{RecolFormel}
Le foncteur $q^* : ({\mathcal Coh}_\X) \rightarrow (\Delta {\mathcal Coh}_\X)$ est une \'equivalence de cat\'egories. Soit $G$ un groupe fini, on a le m\^eme r\'esultat en rempla\c cant la cat\'egorie des modules coh\'e\-rents par celle des $G$-alg\`ebres coh\'erentes galoisiennes (c'est \`a dire des $\O_\X$-alg\`ebre ${\mathcal A}$ munie d'une action de $G$ et telle que ${\mathcal A}^G=\O_\X$).
\end{Pro}

\noindent{\it Preuve :}  Voir \cite{Pri}, Theorem 3.4. \Eproof

\subsection{Rel\`evement des points isol\'es}

On suppose ici ue $R$ est un anneau de valuation discr\`ete complet tel que $\car \ K \ge 0$. Consid\'erons un carr\'e cocart\'esien 
\begin{equation}
\label{CarreCocartesien}
\xymatrix {{k[[x, y]] \over (xy)} \ar[r] &  k((x)) \times k((y)) \\
{k[[x' y']]\over (x'y')} \ar[u]^{\bar\phi} \ar[r] & k((x'))\times k((y')) \ar[u]_{\bar\phi_o \times \bar\phi_t} }
\end{equation}
o\`u $\bar\phi$, $\bar\phi_o$ et $\bar\phi_t$ sont des morphismes de degr\'e $p$. On suppose de plus que $\Spec(\bar\phi_o)$ (resp. $\Spec(\bar\phi_t)$) est un torseur sous un sch\'ema en groupe $G_o$ (resp. $G_t$) de degr\'e $p$, isomorphe \`a $\ZZ/p\ZZ$ ou \`a $\aalpha_p$. Il lui correspond donc un conducteur $m_o$ (resp. $m_t$) premier \`a $p$ (non nul car $G_o$ et $G_t$ ne sont pas isomorphes \`a $\mmu_p$) et un r\'esidu $h_o$ (resp. $h_t$). On dira qu'un tel carr\'e cocart\'esien est {\it $p$-admissible} si on a de plus $m_o=-m_t$, $h_o=-h_t$.

\medskip

Soit $G$ un groupe abstrait $p$-cyclique. \'Etant donn\'e un rev\^etement $p$-admissible d'un point double comme dans \eqref{CarreCocartesien} et un entier $e>0$, une {\it $R$-d\'eformation} d'\'epaisseur $e$ est la donn\'ee d'une action de $G$ sur $R[[x, y]]/(xy-t^e)$, d'une injection $\phi$ de $R[[x', y']]/(x'y'-t^{pe})$ dans $R[[x, y]]/(xy-t^e)$ d'image $\left ( R[[x, y]]/(xy-t^e) \right ) ^G$ et d'un isomorphisme de la fibre sp\'eciale du carr\'e cocart\'esien 
$$\xymatrix{{R[[x, y]] \over (xy-{t}^e)} \ar[r] &  R[[x]]\{x^{-1}\} \times R[[y]]\{y^{-1}\} \\
{R[[x' y']]\over (x'y'-{t}^{pe})} \ar[u]^\phi \ar[r] & R[[x']]\{{x'}^{-1}\}\times R[[y']]\{{y'}^{-1}\} \ar[u]_{\phi_o \times \phi_t} }$$
avec \eqref{CarreCocartesien} (avec les notations de la section pr\'ec\'edente, $\X_a^o\cong Spec \ R[[x]]\{x^{-1}\}$).

\medskip

\begin{Pro}
Soit $K$ un corps complet de caract\'eristique $p > 0$. Soient $m_o > 0$, $e > 0$, $d_o \ge 0$ des entiers, posons $d_t=d_o+em_o(p-1)$, $t$ une uniformisante de $K$ et $\pi=t^{pe}$. Soit $C$ la courbe projective lisse d\'efinie par l'\'equation $V^p-t^{d_t}V+T^{n_t}=0$. Soit ${\mathcal P}$ le mod\`ele
$$\Spec \ R[{1 \over T}] \cup \Spec \ R [{\pi \over T}, T] \cup \Spec \ R[{T \over \pi}].$$
Consid\'erons $f:\C \rightarrow {\mathcal P}$ la normalisation de $\mathcal P$ dans $C$, et $q_0$ le point d'intersection de ${\mathcal P}_s$. Alors $f^{-1}(q_0)$ est un point double ordinaire.
\end{Pro}

\noindent{\it Preuve :} Soient $W_o=f^{-1}(\Spec \ R[{T \over \pi}])$ et $W_t=f^{-1}(\Spec \ R [{1 \over T}, T])$. Alors $W_o \cup W_t$ contient un ouvert dense de $\C_s$. Consid\'erons les $R$-alg\`ebres int\`egres
$A_o=R[{T \over \pi}, U]/(U^p-t^{d_o}U+\left({T \over \pi} \right )^{m_o})$ et $A_t=R[T, {1 \over T}, V]/(V^p-t^{d_t}V+T^{-m_o})$.
 On peut les idetifier \`a des sous alg\`ebres de $K(C)$ en posant $U={V \over \pi^{m_o}}$ et on a alors $\Frac(A_o)= \Frac(A_t)=K(C)$.

De plus, $\Spec \ A_o$ et $\Spec \ A_t$ sont \`a fibres g\'en\'eriques normales et \`a fibre sp\'eciales r\'eduites, ce sont dont des sch\'emas normaux. Comme ils sont finis respectivement sur $\Spec \ R[{T \over \pi}]$ et $\Spec \ R [{1 \over T}, T]$, on a $W_o=\Spec \ A_o$ et $W_t=\Spec \ A_t$. Le sch\'ema $\C$ est de Cohen-Macaulay, donc $\C_s$ v\'erifie la condition (S$_1$) (cf. \cite{LiuBook}, Exercice 8.2.11) et est donc r\'eduite (cf. \cite{LiuBook}, Proposition 7.1.15), r\'eunion de deux composantes irr\'eductibles $\Gamma_o:=\bar W_{o, s}$ et $\Gamma_t:=\bar W_{t, s}$.

On a $f^{-1}(q_0) \subset \Gamma_o \cap \Gamma_t \subset \Gamma_o \setminus W_{o, s}$ et ce dernier est r\'eduit \`a un point (par unicit\'e de la compactification de $W_{o, s}$ par des points r\'eguliers). Par cons\'equent, $f^{-1}(q_0)$ est r\'eduit \`a un point $q$.

Notons $\rho:Z \rightarrow \C_s$ la normalisation en $q$. On a la suite exacte
$$0 \rightarrow \O_{\C_s} \rightarrow \rho_* \O_Z \rightarrow {\mathcal F} \rightarrow 0,$$ 
o\`u $\mathcal F$ est un faisceau gratte-ciel \`a support dans $q$ et avec ${\mathcal F}(\C_s)=\O'_{\C_s, q} / \O_{C_s, q}$ (o\`u $'$ signifie cl\^oture int\'egrale). On en d\'eduit une suite exacte longue
$$0 \rightarrow k \rightarrow k^2 \rightarrow {\mathcal F}(\C_s) \rightarrow \H^1(\C_s, \O_{C_s}) \rightarrow \H^1(Z, \O_Z)\rightarrow 0$$ avec $\dim_k \H^1(\C_s, \O_{C_s})=g(C)$. Calculons $\H^1(Z, \O_Z)$. Soit $\Gamma_o'$ la normalisation de $\Gamma_o$ en $q$, c'est aussi le compl\'et\'e de $W_{o, s}$ par adjonction d'un point r\'egulier. Il est facile de voir que $p_a(\Gamma'_o)=(m_o-1)(p-1)/2=g(C)$.
Or $\H^1(Z, \O_Z) \rightarrow \H^1(\Gamma'_o, \O_{\Gamma'_o})$ est surjectif car $\Gamma'_o$ est une composante connexe de $Z$, il suite que $\dim_k \H^1(Z, \O_Z) \ge g(C)$, et donc $\dim_k {\mathcal F}(\C_s)=1$. Cela implique que $q$ est un point double ordinaire. \Eproof

\begin{Cor}
Soient $d_o\ge 0$, $d_t \ge 0$ et $e>0$ des entiers tels que $d_t-d_o=em_o(p-1)$. On suppose que $d_o$ (resp. $d_t$) est nul si et seulement si $G_o$ (resp. $G_t$) est \'etale. Alors il existe une $R$-d\'eformation du carr\'e cocart\'esien \eqref{CarreCocartesien} telle que les diff\'erentes aux bords soient $d_o$ et $d_t$.
\label{Relevement_point_double}
\end{Cor}

\noindent{\it Preuve :}  
Il suffit de consid\'erer le morphisme $\hat \O_{{\mathcal P}, q_0} \rightarrow \hat \O_{\C, q}$ de la proposition pr\'ec\'edente. \Eproof
\medskip

Int\'eressons nous maintenant au relevement des points r\'eguliers ramifi\'es.

\begin{Pro}
\label{Relevement_Point_Lisse}
Soit $f:\Spec k[[u]] \rightarrow \Spec k[[v]]$ un morphisme tel que le morphisme induit $\Spec k((u)) \rightarrow \Spec k((v))$ soit un $\aalpha_p$-torseur d\'efinit par une fonction $w \in k((v))$ de valuation $1-pm$ avec $m \ge 0$. Soit $n >0$ un entier. Alors il existe un morphisme $\Spec R[[U]] \rightarrow \Spec R[[V]]$ tel que le morphisme induit $\Spec R[[U]]\{1/U\} \rightarrow \Spec R[[V]]\{1/V\}$ soit un $\M_n$-torseur et dont la fibre sp\'eciale est $f$. 
\end{Pro}

\noindent{\it Preuve :} Par d\'efinition de $w$, $k((u))$ est de la forme $k((v))[u]/(u^p-w)$. Quitte \`a changer de param\`etre $v$, on peut supposer que $w=v^{1-pm}$. Par suite, $$R[[V]][U]/(U^p-t^{(p-1)n}V^{(p-1)m}U+V)$$ convient. \Eproof

\subsection{Preuve du th\'eor\`eme \ref{Relevement}}

Cette d\'emonstration s'inspire de celle du th\'eor\`eme de rel\`evemnt de \cite{HenrioNP}.

Par d\'efinition d'un graphe de Hurwitz adapt\'e, on a une inclusion canonique $\Gamma_{\Y_s} \rightarrow \Gamma$, on considerera donc les sommets de $\Gamma_{\Y_s}$ \'egalement comme des sommets de $\Gamma$. On note $y_1, \ldots, y_n$ les points lisses de $\Y_s$ en lesquels $f_s$ n'est pas un torseur et $x_i=f_s(y_i)$. On reprend alors les notations de la proposition \ref{RecolFormel}.

Comme $f_s$ n'est pas un hom\'eomorphisme, il existe un sommet $\tilde v \in \Gamma_{\Y_s}$ tel que la restriction de $f_s$ \`a un ouvert non vide de $F_{\tilde v}^0$ soit \'etale. Par suite, $d(\tilde v) = 0$. Les conditions de compatibilit\'es sur les diff\'erentes imposent alors que pour tout sommet $\tilde v_1$ de $\Gamma$, on ait $(p-1) | d(\tilde v_1)$.
Soit $v$ un sommet de $\Gamma_{\X_s}$ et $\tilde v$ un sommet de $\Gamma_{\Y_s}$ au dessus de $v$. Notons $n_v=d(\tilde v)/(p-1)$. Cet entier est ind\'ependant du choix du rel\`evement $\tilde v$ car $\Gamma$ est un graphe de Hurwitz. Supposons dans un premier temps que $\tilde v$ est fixe sous l'action de $G$. D'apr\`es le corollaire \ref{ReleveTorseur}, il existe un $\M^{n_v}$-torseur $U_v$ au dessus de $\X_v$ relevant $f_s$ au dessus de $\X_{v, s}$. Supposons maintenant que $\tilde v$ n'est pas fixe sous l'action de $G$. Posons $U_v:=G \times \X_v$. On a alors une projection canonique $U_v \rightarrow \X_v$. Dans tout les cas, $U_v \rightarrow \X_v$ est un $G$-rev\^etement galoisien. 

Soit $i\in\{1, \ldots, n\}$. Alors, par l'hypoth\`ese (H$_2$), $f_s$ est un hom\'eomorphisme au voisinage de $y_i$. Notons $v_i$ le sommet correspondant \`a la composante contenant $x_i$ et $\tilde v_i$ le sommet de $\Gamma$ au dessus de $v_i$. L'hypoth\`ese (H$_2$) et la proposition \ref{Relevement_Point_Lisse} permettent de relever $f_{s, y}:\Spec \hat\O_{\Y_s, y} \rightarrow \Spec \hat\O_{\X_s, x_i}$ et ainsi d\'efinir un $G$-rev\^etement galoisien $U_i$ de $\X_i$ telle que la diff\'erente en le point g\'en\'erique de $\X_{i, s}$ soient $d(\tilde v)=(p-1)n_i$. Construisons maintenant l'isomorphisme $\phi_i$. On sait que 
$$\O_{U_{v_i}, \eta_i} \otimes_{\O_{\X, \eta_i}} \widehat{(\O_{\X, x_i})_{\eta_i}}\cong R[[T]]\{T^{-1}\}[U]/(U^p-t^{n_i}U-Q)$$ avec $Q \in R[[T]]\{T^{-1}\}$ et 
$$\O_{U_i, \eta_i}\otimes_{\O_{\X, \eta_i}} \widehat{(\O_{\X, x_i})_{\eta_i}}
\cong R[[T]]\{T^{-1}\}[U']/({U'}^p-t^{n_i}U-Q')$$ avec $Q' \in R[[T]]\{T^{-1}\}$ et $Q=Q'\mod t$. Par suite, $Q=T^{-m}(\alpha+tQ_1(1/T)+TQ_2(T)+tQ_3(T))$ avec $Q_1\in R\{T\}$, $Q_2, \ Q_3 \in R[[T]]$, $\alpha \in R \setminus tR$ et $m >0$ un entier premier \`a $p$. Posons $T'=T(\alpha+tQ_1(1/T)+TQ_2(T)+tQ_3(T))^{1 \over m}$. Alors $T'$ est un param\`etre local de $\X$ en $x_i$ et on a donc 
$$\O_{U_{v_i}, \eta_i} \otimes_{\O_{\X, \eta_i}} \widehat{(\O_{\X, x_i})_{\eta_i}}\cong R[[T']]\{{T'}^{-1}\}[U]/(U^p-t^{n_i}U-{T'}^{-m}).$$ De m\^eme, on peut trouver un param\`etre $T''$ de $R[[T]]\{T^{-1}\}$ tel que
$$\O_{U_i, \eta_i}\otimes_{\O_{\X, \eta_i}} \widehat{(\O_{\X, x_i})_{\eta_i}}
\cong R[[T'']]\{{T''}^{-1}\}[U']/({U'}^p-t^{n_i}{U'}-{T''}^{-m}).$$ Le morphisme $\phi_i$ est alors d\'efini par $\phi_i(T')=T''$ et $\phi_i(U)=U'$.

Soit $a$ une ar\^ete de $\Gamma_{\X_s}$ et $\tilde a$ une ar\^ete de $\Gamma_{\Y_s}$ au-dessus de $a$. Supposons dans un premier temps que $\tilde a$ est fixe sous l'action de $G$, la proposition \ref{Relevement_point_double} permet de construire un rel\`evement $U_a$ de $f_s$ au dessus de $\X_a$ avec, au bord de la couronne formelle, les diff\'erentes $d(o(\tilde a))$ et $d(t(\tilde a))$. Par suite, on peut contruire des isomorphismes $\phi_a$ et $\phi_{\bar a}$ par la m\^eme m\'ethode que pour les points $x_i$.
Supposons maintenant que $\tilde a$ n'est pas fixe sous l'action de $G$. On pose $U_a:=G \times \X_a$. Supposons que $o(a)$ soit non fixe sous l'action de $G$, on peut alors d\'efinir l'isomorphisme $\phi_a$ comme provenant de l'isomorphisme $\alpha_a^* \O_{\X_a} \rightarrow \beta_a^* \O_{\X_{o_a}}$.
Supposons maintenant que $o(a)$ est fixe sous l'action de $G$. Comme $a$ n'est pas fixe, l'action de $G$ est non triviale sur $o(a)$ et on a un isomorphisme de $\O_{\X_{o(a)}}[G]-$module $\O_{U_{o(a)}} \otimes _{\O_{\X_{o(a)}}} \O_{\X_a^o} \cong G \times \O_{\X_a^o}$. On obtient donc un isomorphisme, not\'e $\phi_a$,
$$\O_{U_{o(a)}} \otimes _{\O_{\X_{o(a)}}} \O_{\X_a^o} \cong G \times \O_{\X_a^o} \cong G \times \O_{\X_a} \otimes_{\O_{\X_a}} \O_{\X_a^o} = U_a \otimes_{\O_{X_a}} \O_{\X_a^o}.$$

Par suite, le th\'eor\`eme \ref{RecolFormel} nous permet de construire une $G$-alg\`ebre galoisienne $\mathcal A$ de sorte que $\Spf(\mathcal A)$ r\'eponde \`a la question.
\Eproof

\begin{Rem}
On peut remarquer que dans les donn\'ees du th\'eor\`eme pr\'ec\'edent, il n'y a pas de condition de cocycle \`a v\'erifier.
\end{Rem}

\section{Un th\'eor\`eme de rel\`evement des courbes hyperelliptiques}

On fixe un anneau de valuation discr\`ete complet $R$. On souhaite montrer dans cette section que toute courbe  stable hyperelliptique $C$ sur $\Spec(k)$ (cf. d\'efinition \ref{CourbeStableHyperelliptique}) s'\'etend, quitte \`a faire une extension finie de $R$, en une courbe hyperelliptique $\C$ sur $\Spec(R)$ dont la fibre g\'en\'erique est lisse, et ceci ind\'ependamment des caract\'eristiques de $K$ et $k$. Les techniques utilis\'ees n'\'etant pas diff\'erentes de celles utiles pour le rel\`evement des rev\^etements $p$-admissibles en g\'en\'eral, on traitera ce dernier cas.

Le cas o\`u $\car \ k \not = p$ est trait\'e dans \cite{Hen_galoisien} Th\'ero\`eme 3.7, on supposera d\'esormais que $\car \ k = p$. On \'etudira d'abord le cas d'\'egale caract\'eristique et on montrera que cela implique le cas d'in\'egale caract\'eristique.

\subsection{R\'eduction des hypoth\`eses du th\'eor\`eme \ref{Relevement}}

On suppose ici que $\car k=p$. Comme on ne peut pas d\'eformer $\mmu_p$ en $\ZZ/p\ZZ$ en \'egale caract\'eristique $p$, il est n\'ecessaire de se restreindre au cas o\`u les torseurs intervenant dans les rev\^etements $p$-admissibles sont des $\aalpha_p$ ou des $\ZZ/p\ZZ$-torseurs. 

\begin{The}
\label{ChangeTorseur}
Soit $k$ un corps alg\'ebriquement clos de caract\'eritique $p>0$. Soit $(f:Y \rightarrow X, \{(G_v, u_v)\})$ un rev\^etement $p$-admissible de courbes semi-stables sur $k$. Supposons que $p_a(X)=0$. Alors il existe un rev\^etement $p$-admissible $(f:Y \rightarrow X, \{(G'_v, u'_v)\})$ tel que pour tout $v$, on a $G'_v\not = \mmu_p$ (i.e. on a \'elimin\'e les $\mmu_p$-torseurs du rev\^etement $p$-admissible).

Plus pr\'ecis\'ement, si $G_v\not = \mmu_p$ alors $(G'_v, u'_v)=(G_v, u_v)$. Si $G_v=\mmu_p$, on a $G'_v=\aalpha_p$, et pour toute ar\^ete $a$ d'origine $v$ telle que $m(a)\not = 0$, le conducteur $m'(a)$ associ\'e \`a $G'_v$ en $a$ v\'erifie $m'(a)=m(a)$.
\end{The}

\noindent{\it Preuve :} D'apr\`es le lemme \ref{NN} (c), quite \`a remplacer notre rev\^etement $p$-admissible par un autre qui lui est \'equivalent, on peut supposer que pour tout sommet $v$ avec $G_v=\mmu_p$, et pour tout point $P$ de $F_v$ tel que $m(P) \not = 0$, on a $\Ord_{f(P)} u_v=0$.

On peut supposer qu'il existe un sommet $v_0$ tel que $G_{v_0}=\mmu_p$. Soit $d(,)$ la distance sur le graphe $\Gamma_Y$ (elle existe m\^eme si $\Gamma_Y$ n'est pas forc\'ement un arbre). Nous allons construire les $u'_v$ par r\'ecurrence sur $d(v_0, v)$, $G'_v$ \'etant toujours d\'efinit comme dans l'\'enonc\'e du th\'eor\`eme. Posons $u'_{v_0}=u_{v_0}$. Pour toute ar\^ete $a$ d'origine $v_0$ telle que $m(a) \not = 0$, on a $\Ord_{a} du'_{v_0}=\Ord_a({du_{v_0} \over u_{v_0}})$ donc $m'(a)=m(a)$. Fixons un entier $d \ge 0$. Supposons avoir construit $u'_v$ pour tout $v$ de distance $d(v_0, v) \le d$, avec les propri\'et\'es suivantes :
\begin{enumerate}[(A)]
\item $u'_v=u_v$ si $G \not = \mmu_p$;
\item si $G_v=\mmu_p$, alors $u'_v/u_v$ est une puissance $p$-i\`eme, et si $a$ est une ar\^ete d'origine $v$ tel que $m(a) \not = 0$, alors $m'(a)=m(a)$;
\item la condition de compatibilit\'e $m'(a)=-m'(\bar a)$ est v\'erifi\'ee pour toute ar\^ete ramifi\'ee reliant les sommets $v$ de distance $d(v, v_0)\le d$ (noter que $h'(a)=h'(\bar a)=0$ car $G'_v=\aalpha_p$).
\end{enumerate}

Soit $v$ un sommet de $\Gamma_Y$ tel que $d(v, v_0) = d+1$. Supposons que $F_v$ rencontre les composantes $F_{\tilde v}$ avec $d(\tilde v, v_0) \le d+1$ en au moins deux points, alors $F_v$ est contenu dans un cycle de $\Gamma_Y$. Comme $\Gamma_X$ est un arbre, il suit que les points d'intersection de $F_v$ avec ces autres composantes s'envoient sur le m\^eme point dans $X$, et sont donc non ramifi\'es pour $f$. Par suite, $v$ ainsi que les sommets $\tilde v$ tels que $d(\tilde v, v_0)\le d+1$ et que $F_v \cap F_{\tilde v} \not = \emptyset$ sont non ramifi\'es, et il n'y a pas d'ar\^ete ramifi\'ee reliant $v$ et $\tilde v$. On prend donc $u'_v=u_v$. Les condictions (A)-(C) ci-dessus sont triviales.

On peut donc supposer que $v$ est reli\'e \`a un unique sommet $\tilde v$ tel que $d(\tilde v, v_0) \le d+1$, et que $F_v \cap F_{\tilde v}=\{P\}$. On a alors $d(\tilde v, v_0)=d$. Soit $a$ l'ar\^ete reliant $v$ \`a $\tilde v$. Supposons d'abord que $m(a)\not = 0$. Posons $u'_v=u_v$. Alors $m'(a)=m(a)$ (comme pour $v_0$ si $G_v=\mmu_p$). Par hypoth\`ese de r\'ecurrence, $m'(\bar a)=m(\bar a)$ car ce dernier est non nul, il suit que $m'(a)=-m'(\bar a)$. On a donc les conditions (A)-(C).

Il reste le cas o\`u $m(a)=0$. On a alors $G_v=\mmu_p$ et $f(F_v)=\PP^1_k$. De m\^eme pour $G_{\tilde v}$. La condition (B) sur $\tilde v$ implique que $\Ord_P u'_v$ est congru \`a $\Ord_P u_{\tilde v}$ modulo $p$. Comme $h(a)=-h(\bar a)$, il suit du lemme \ref{NN} (b) que $n'_{\tilde v}:=\Ord_P u'_{\tilde v}$ est congru \`a $-n_v$ modulo $p$, o\`u $n_v:=\Ord_P u_v$. Il existe $Q \in f(F_v)$ tel que $\Ord_Q u_v$ soit premier \`a $p$ (car $n_v$ est premier \`a $p$). Choisissons un param\^etre $x$ de $f(F_v)$ avec $(x)=P-Q$. Soient $r=(n'_{\tilde v}+n_v)/p \in \ZZ$ et $u'_v=x^{pr}u_v$. alors $m'(\bar a)=-n'_{\tilde v}=1+\Ord_P du'_v=-m'(a)$. On a donc v\'erifi\'e la condition (C). Soit $a_1$ une ar\^ete d'origine $v$ avec $m(a_1)\not = 0$. Alors le point $P_1$ correspondant \`a $a_1$ est distinct de $P$ et $Q$. Il suit imm\'ediatement du fait $x, u_v \in \O^*_{f(P_1)}$ que 
$$m'(a_1)=-(1+\Ord_{P_1} du'_v)=-(1+\Ord_{P_1} du_v)=m(a_1).$$
Ce qui ach\`eve la v\'erification de la condition (B) ainsi que la preuve de ce th\'eor\`eme. \Eproof

\begin{Pro}
\label{ChangeTorseur2}
Soit $(f:Y \rightarrow X, \{(G_v, u_v\}_v)$ un rev\^etement $p$-admissible tel que pour tout sommet $v$ de $\Gamma_Y$, $G_v$ soit isomorphe \`a $\ZZ/p\ZZ$ o\`u $\aalpha_p$. Alors il existe un rev\^etement $p$-admissible $(f':Y' \rightarrow X', \{(G'_v, u'_v\}_v)$ v\'erifiant les hypoth\`ese {\rm (H$_1$)} et {\rm (H$_2$)} du th\'eor\`eme \ref{Relevement} et des immersions ferm\'ees $Y \rightarrow Y'$ et $X \rightarrow Y'$ avec $p_a(X')=p_a(X)$, $p_a(Y')=p_a(Y)$. 
\end{Pro}

\noindent{\it Preuve :} Notons $Y'$ (resp. $X'$) la courbe semi-stable obtenue en ajoutant transversalement \`a $Y$ (resp. $X$) une composante $\PP^1_k$ en chaque point r\'egulier $P$ en lequel $f$ n'est pas un torseur (resp. en $f(P)$). Pour chaque sommet $v'$ de $\Gamma_{Y'}$ provenant d'un sommet $v$ de $\Gamma_Y$, posons $(G'_{v'}, u'_{v'}):=(G_v, u_v)$. Soit $v'$ un sommet de $\Gamma_{Y'}$ ne provenant pas d'un sommet de $\Gamma_Y$. Posons $G'_{v'}=\aalpha_p$. Par construction, il existe un unique sommet $\tilde v$ de $\Gamma_{Y'}$ tel que $F_{v'} \cap F_{\tilde v} \not = \emptyset$ et cette intersection est r\'eduite \`a un point $P$. Notons $n:=\Ord_P u_{\tilde v}$. On montre ais\'ement qu'il existe un entier $q > 0$ et des entiers $n_1, \ldots, n_q \ge 0$ tels que $$n=\sum_{i=1}^q 1-n_i p.$$ Notons $V$ la composante ajout\'ee  \`a $X$ en $f(P)$. Soient $P_1, \ldots, P_q$ des points de $V$ distincts deux \`a deux et de $f(P)$. Alors il existe $u_{v'} \in \H^0(V\setminus\{f(P), P_1, \ldots, P_q\}, \O_V^*)$ tel que $(u_{v'})=-nf(P)+\sum_i (1-pn_i)P_i$. Cette section rationnelle permet alors de munir un ouvert dense de $F_{v'}^0$ d'une structure de $G_{v'}$ torseur telle que l'hypoth\`ese (H$_2)$ soit v\'erifi\'ee et ceci permet de d\'efinir un morphisme $f':Y' \rightarrow X'$. On montre alors ais\'ement que $(f':Y' \rightarrow X', \{(G'_{v'}, u'_{v'})\}_{v'})$ est un rev\^etement $p$-admissible. \Eproof

\subsection{Construction des donn\'ees de Hurwitz}

\begin{Pro}
\label{DonHurwitz}
Soient $(f:Y \rightarrow X, \{(G_v, u_v)\}_v)$ un rev\^etement $p$-admissible avec $p_a(Y) \ge 2$ et $p_a(X)=0$. On suppose que si la caract\'eristique de $K$ est $p$, la description du rev\^etement $p$-admissible $f$ ne fait pas intervenir de $\mmu_p$-torseurs. Alors, quitte \`a faire une extension finie de $K$, il existe un mod\`ele $\X$ de $X$ sur $R$ et un graphe de Hurwitz $\Gamma$ adapt\'e \`a $\X$ et \`a $(f:Y \rightarrow X, \{(G_v, u_v)\}_v)$.
\end{Pro}

\noindent{\it Preuve :} 
Notons $\Gamma$ le graphe orient\'e (l'orientation \'etant choisie arbitrairement) dont l'ensemble des sommets est la r\'eunion de l'ensemble des sommets de $\Gamma_Y$ et de l'ensemble des points de $Y$ en lesquels $f$ n'est pas un torseur (cf. Exemple \ref{ExempleReduction}). Le graphe $\Gamma$ vient alors naturellement avec des applications $m$, $h$ v\'erifiant les conditions de compatibilit\'es (cf. d\'efinition \ref{DonneeCompatible}). Soit $v$ un sommet de $\Gamma$, si $v$ correspond \`a une composante irr\'eductible $Y_1$ de $Y$, on pose $g(v)=p_a(Y_1)$ et sinon on pose $g(v)=0$.

Par suite, les diff\'erentes \'egalit\'es ne faisant intervenir que les fonctions $g$, $m$ et $h$ proviennent des \'egalit\'es classiques mentionn\'ee dans la d\'efinition des graphes de Hurwitz.

Il reste alors \`a d\'efinir les fonction $\varepsilon$ et $d$. 
\'Etant donn\'e que la courbe $Y$ est stable de genre au moins $2$, il y a au moins un sommet \'etale $v_0$, \'eventuellement non fixe sous $G$ (sinon $Y$ serait de genre nulle). On peut choisir un ordre partiel (comme dans la d\'emonstration de \ref{ChangeTorseur}) sur $\Gamma'$ (car celui-ci est un arbre) tel que $v_0$ soit le plus petit \'el\'ement pour cette ordre. Nous allons construire les fonctions $\varepsilon$ et $d$ par ``r\'ecurrence'' sur cet ordre partiel.

Posons $d(v_0)=0$. 
Soit $a$ une ar\^ete de $\Gamma_X$ joignant $v_1$ \`a $v_2$ avec $v_1 \le v_2$. Supposons que $d$ soit d\'ej\`a d\'efinit en $v_1$, nous allons la d\'efinir en $v_2$.

\begin{enumerate}[-]
\item Si $v_2$ est \'etale, posons $d(v_2)=0$. Comme on a  $m(\bar a) > 0$ (car $\bar a$ est une ar\^ete d'origine \'etale) et que $d(v_1) \ge 0$, on peut prendre $\varepsilon(a) =-{d(v_1) \over m(a)(p-1)} \ge 0$.
\item Si $v_2$ est additif
\begin{enumerate}[-]
\item Si $v_1$ est \'etale alors $m(a) > 0$ et $d(v_1)=0$, on peut donc trouver  $\varepsilon(a)$ tel que $0<(p-1)\varepsilon(a) m(a)$ et par suite, on pose $d(v_2)=(p-1)\varepsilon(a) m(a)$
\item Si $v_1$ est additif on a $0 < d(v_1) $ donc on peut prendre $\varepsilon(a)$ tel que $0 < d(v_1)+(p-1)m(a)\varepsilon(a) $, on pose alors $d(v_2)=d(v_1)+m(a)\varepsilon(a)(p-1)$
\end{enumerate}
\end{enumerate}

Dans tout les cas, on peut trouver $d(v_2)$ et $\varepsilon(a)$ tel que 
$$d(v_2)-d(v_1)=(p-1)m(a) \varepsilon(a).$$ On peut donc construire par r\'ecurrence des applications $\varepsilon$ et $d$ \`a valeur rationnelle v\'erifiant les hypoth\`eses. Par suite, quitte \`a multiplier $\varepsilon$ et $d$ par une m\^eme constante, on peut supposer qu'elles sont \`a valeures enti\`eres. En relevant ces applications sur $\Gamma$, on fait de celui-ci un graphe de Hurwitz. 

On peut alors construire une mod\`ele $\X$ sur $R$ de $X$ tel $\Gamma$ soit un graphe de Hurwitz adapt\'e \`a $\X$ et \`a $(f, \{(G_v, u_v)\}_v)$ d'apr\`es \cite{Hen_galoisien} Proposition 2.6. \Eproof

\subsection{Th\'eor\`eme de rel\`evement des courbes hyperelliptiques stables}
\begin{The}
\label{Releve}
Soit $R$ un anneau de valuation discr\`ete complet d'\'egale caract\'e\-ristique $p$.
Soit $C \rightarrow \Spec(k)$ une courbe stable hyperelliptique. Alors, quitte \`a faire une extension finie de $R$, il existe une courbe $\C \rightarrow \Spec(R)$ hyperelliptique semi-stable dont la fibre g\'en\'erique est lisse et dont la fibre sp\'eciale est isomorphe \`a $C$.
\end{The}

\noindent{\it Preuve :} 
Dans le cas o\`u $C$ est lisse il suffit de relever une \'equation de $k$ \`a $R$. On peut donc supposer que $C$ n'est pas lisse.

D'apr\`es la proposition \ref{Reduction_Admissible}, la courbe $C$ vient naturellement avec un rev\^etement $2$-admissible $f:Y \rightarrow X$ si $\car\  k \not = 2$ (ou $(f:Y \rightarrow X, \{(G_v, u_v)\})$ si $\car k = 2$) o\`u $p_a(Y)=p_a(C)$ et $Y$ est obtenue de $C$ en ajoutant des composantes $\PP^1_k$.
 
Supposons tout d'abord que la caract\'eristique de $k$ est 2. D'apr\`es la r\'eduction pr\'ec\'edente (th\'eor\`emes \ref{ChangeTorseur} et \ref{ChangeTorseur2}), on peut supposer que le rev\^etement $(f, \{(G_v, u_v)\}_v)$ v\'erifie les hypoth\`eses (H$_1$) et (H$_2$) du th\'eor\`eme \ref{Relevement}. Par suite, gr\^ace \`a la proposition \ref{DonHurwitz}, on peut construire, quitte \`a faire une extension finie de $K$, un mod\`ele $\X$ de $X$ et un arbre de Hurwitz $\Gamma$ adapt\'ee au rev\^etement \`a $\X$ et \`a $(f, \{(G_v, u_v)\}_v)$. On applique alors le corollaire \ref{RelevementCourbe} qui nous fournit le rel\^evement $\Y \rightarrow \X$. Le mod\`ele stable $\C$ de $\Y_K$ r\'epond \`a la question.

Regardons maintenant le cas o\`u la caract\'eristique de $k$ est diff\'erente de $2$. On a une action kummerienne de $G=\ZZ/2\ZZ$ sur $Y$ dont le quotient est $X$. Comme le cardinal de $G$ est premier, on peut appliquer le th\'eor\`eme 3.7 de \cite{Hen_galoisien} qui nous prouve que cette action se rel\`eve en un mod\`ele semi-stable $\Y \rightarrow \Spec(R)$. Finalement, le mod\`ele stable $\C$ de $\Y_K$ r\'epond \`a la question.  \Eproof

\begin{Cor}
\label{ReleveCarZero}
Soient $k$ un corps alg\'ebriquement clos et $C \rightarrow k$ une courbe stable hyperelliptique. Alors il existe un anneau de valuation discr\`ete $R$ d'in\'egales caract\'eristiques  et dont le corps r\'esiduel est $k$, et une courbe $\C \rightarrow \Spec \ R$ dont la fibre g\'en\'erique est lisse et hyperelliptique et dont la fibre sp\'eciale s'identifie \`a $C$.
\end{Cor}

\noindent{\it Preuve :} D'apr\`es le th\'eor\`eme pr\'ec\'edent, on peut relever $C$ en une courbe $\C_1 \rightarrow \Spec \ R_1$ o\`u $R_1$ est un anneau de valuation discr\`ete d'\'egale caract\'eristique. D'autre part, on peut relever la fibre g\'en\'erique de $\C_1$ en une courbe lisse sur un anneau de valuation discr\`ete $R_2$ dont le corps des fractions est de caract\'eristique z\'ero (car il suffit de relever une \'equation hyperelliptique). Ceci prouve que le point g\'en\'erique de l'anneau local du point correspondant \`a $C \rightarrow \Spec \ k$ dans $\bar I_g$ correspond \`a une courbe lisse en caract\'eristique z\'ero. Le r\'esultat est alors direct. \Eproof

\begin{Rem}
Le th\'eor\`eme \ref{Releve} et le corollaire \ref{ReleveCarZero} sont encore vrais pour des rev\^etements $p$-admissibles d'une courbe de genre z\'ero, la preuve du corollaire \ref{ReleveCarZero} devant \^etre adapt\'ee  car il n'y a pas, a priori, d'espace de module jouant le role de $\bar I_g$.
\end{Rem}

\begin{The}
Soit $S$ le spectre d'une anneau de valuation discr\`ete $R$.
On a l'\'egalit\'e (ensembliste) $\bar I_{g, S} \times_ S \Spec \ \kappa(s) = \bar I_{g, s}.$ En particulier,  $\bar I_{g, S} \times_ S \Spec\  \kappa(s)$ est irr\'eductible.
\end{The}

\noindent{\it Preuve :} Quitte \`a faire un changement de base plat, on peut supposer $R$ complet. La premi\`ere partie du th\'eor\`eme d\'ecoule directement de ce qui pr\'ec\`ede. En effet, par construction on a une immersion ferm\'ee $\bar I_{g, s} \rightarrow \bar I_{g, S} \times_S \Spec \ \kappa(s) $. Soit $\p$ un point de $\bar I_{g, S} \times_S \Spec \ \kappa(s)$. Alors $\p$ correspond \`a une courbe stable hyperelliptique (car sp\'ecialisation d'un point correspondant \`a une courbe lisse hyperelliptique). Le th\'eor\`eme \ref{Releve} implique alors que ce point est dans l'adh\'erence du lieu des courbes lisses hyperelliptiques sur $\kappa(s)$. C'est donc un point de $\bar I_{g, s}$.

 L'irr\'eductibilit\'e de $\bar I_{g, S} \times_ S \Spec \ \kappa(s)$ provient alors de l'irr\'eductibilit\'e de $I_{g, \Spec \ \kappa(s)}$ montr\'ee dans \cite{Lo} Theorem 4.1. \Eproof

\appendix \section{Calcul en caract\'eristique diff\'erente de 2}

Dans toute cette section, $R$ sera un anneau de valuation discr\`ete de caract\'eristique r\'esiduelle diff\'erente de 2. On note $S=\Spec(R)$, $\C \rightarrow S$ une courbe stable hyperelliptique de genre $g \ge 2$ dont la fibre g\'en\'erique est lisse et dont la fibre sp\'eciale poss\`ede deux composantes irr\'eductibles lisses ou bien une seule composante irr\'eductible avec un point singulier. Sous ces hypoth\`eses, on se propose de d\'emontrer explicitement les in\'egalit\'es (\ref{TheoPrinc}). 
Ce calcul reprend les id\'ees expos\'ees dans \cite{Kau}. Toutefois, la partie combinatoire n'est pas utile ici car les courbes sont suffisament simples.

\begin{Pro}
Si $\C_s$ poss\`ede un unique point singulier de type $\beta_j$ (cf. la discussion suivant la {\rm remarque} $5$) alors $$\Ord_s(\Lambda)=4(g-j)j\delta_s.$$
Si $\C_s$ poss\`ede deux points singuliers de type $\alpha_j$ alors $$\Ord_s(\Lambda)=(g-j)(j+1)\delta_s.$$
\end{Pro}

\noindent{\it Preuve :}  Soit $\sigma$ l'involution hyperelliptique de $\C_\eta$. Par unicit\'e du mod\`ele stable, $\sigma$ se prolonge en une involution de $\C$. On peut remarquer que dans le cas qui nous int\'eresse, les composantes irr\'eductibles de $\C_s$ sont fixes par $\sigma$. Notons $Y=\C/\gsigma$. Quitte \`a faire une extension finie du corps de base $K$, on peut supposer que les points fixes de $\C_\eta$ sous l'involution hyperelliptique sont rationnels, notons les $a_1, \ldots, a_{2g+2}$. Par suite, notons $b_1, \ldots ,b_{2g+2}$ leurs images dans $Y_\eta$. 

\begin{Lem}
\label{SpecialisationSimple}
Les sp\'ecialisations des $a_i$ qui sont dans le lieu lisse de $\C_s$ sont deux \`a deux distinctes. 
\end{Lem}

\noindent{\it Preuve :} En effet, au voisinage d'un point $a_i$ qui se sp\'ecialise dans le lieu lisse, $\C$ est de la forme $\Spec R[x, y]/(y^2-Q(x))$, $b_i$ correspondant \`a un z\'ero de $Q$. Si $b_i$ et $b_j$ ont m\^eme sp\'ecialisation $\q$, la courbe $\C$ sera non singuliere en $\q$ si et seulement si $i=j$. \Epl

\begin{Lem}
Si les points singuliers ne sont pas de type $\alpha_0$, les $b_i$ se sp\'ecialisent dans le lieu lisse de $Y_s$.
\end{Lem} 

\noindent{\it Preuve :} Si les points singuliers sont du type $\alpha_j$ ($j > 0$), alors ils ne sont pas fixes par $\sigma$ et donc ne sont pas sp\'ecialisations de $b_i$. Supposons donc que le point singulier de $\C_s$ est du type $\beta_j$. Raisonnons par l'absurde et supposons que $b_i$ se sp\'ecialise dans un point singulier. Par suite, le point $a_i$ se sp\'ecialise en un point singulier $\p$ de $\C_s$. Notons $\tilde\C$ l'\'eclat\'e de $\C$ en $p$ de sorte que $a_i$ se sp\'ecialise en un point lisse de $\tilde\C_s$. Comme $\p$ est fixe par l'involution hyperelliptique, $\sigma$ se prolonge en une involution de $\tilde \C$. Notons $\tilde Y=\tilde \C / \gsigma$. Ainsi d\'efinit, $\tilde Y$ est obtenu de $Y$ en \'eclatant l'image de $\p$ dans $Y_s$. Notons $E$ la transform\'ee stricte de $\p$ par le morphisme $\tilde \C \rightarrow \C$ et $F$ son image dans $\tilde Y$. Le morphisme $\tilde \C \rightarrow \tilde Y$ d\'efinit alors un morphisme g\'en\'eriquement s\'eparable de degr\'e deux $E \rightarrow F$ (car la caract\'eristique r\'esiduelle de $R$ est diff\'erente de 2 et $\sigma$ est non triviale sur $E$). Notons $\alpha$ et $\beta$ l'intersection de $E$ avec le lieu singulier de $\tilde C_s$. Les points $\alpha$, $\beta$ et la sp\'ecialisation de $b_i$ sur $E$ sont alors fixes par l'involution $\sigma$. Le morphisme $E \rightarrow F$ poss\`ede alors au moins trois points de ramification, ce qui contredit le fait que $E$ est de genre nulle. \Eproof

\begin{Lem}
Si le point singulier est de type $\alpha_0$ alors il y a exactement deux $b_i$ se sp\'ecialisant dedans.
\end{Lem}

\noindent{\it Preuve :}  Pour des questions de genre de la composante irr\'eductible de $\C_s$ (formule de Hurwitz), au moins deux $b_i$ se sp\'ecialisent dans le point singulier. Si on suppose qu'au moins trois $b_i$ se sp\'ecialisent dans le lieu singulier de $Y_s$, le m\^eme raisonnement que dans la d\'emonstration du lemme pr\'ec\'edent conduit \`a une contradiction sur le genre de $\C_s$. \Eproof
\medskip

Dans le cas  o\`u le point singulier est de type $\alpha_0$ on peut, grace \`a ce lemme, se ramener \`a un calcul analogue \`a celui des courbes ayant des points singuliers du type $\alpha_j$ avec $j > 0$ en d\'eployant les sp\'ecialisations des points $b_i$ se sp\'ecialisant dans le lieu singulier.

Notons $Y_1$ et $Y_2$ les deux composantes irr\'eductibles de $Y_s$. Si on regarde la courbe $Y$ du point de vue de la g\'eom\'etrie rigide, on peut consid\'erer $Y_1$ comme l'ensemble des points de modules $1$ et $Y_2$ comme l'ensemble des points de modules $|t^b|$ o\`u $b$ est l'\'epaisseur du point singulier de $Y_s$ dans $Y$ et $t$ est une uniformisante de $R$. Ce choix correspond alors \`a un certain param\`etre $x$ de $Y_\eta$ qu'on fixe pour le reste de cette d\'emonstration. 

Notons $y^2+Q(x)y-P(x)$ une \'equation affine de $\C_\eta$ (le param\`etre $x$ \'etant le param\`etre fix\'e precedement). Dire que les $b_i$ sont rationnels \'equivaut \`a dire que $f=Q^2+4P$ est d\'ecompos\'e. Notons $f=A \prod(x-c_i)$. On va alors \'evaluer les diff\'erents facteurs de $\Lambda$. Dans le cas d'une courbe ayant un point du type $\beta_j$, on a une configuration donn\'ee par la figure \ref{PointBetaj} (cf \cite{Co-Ha}, \S 4) et dans le cas d'une courbe ayant deux points de type $\alpha_j$ (\'eventuellement $j=0$ quitte \`a faire un \'eclatement) on a une configuration donn\'ee par la figure \ref{PointAlphaj}.

\begin{figure}
\caption{Point de type $\beta_j$}
\label{PointBetaj}
\centerline{\input{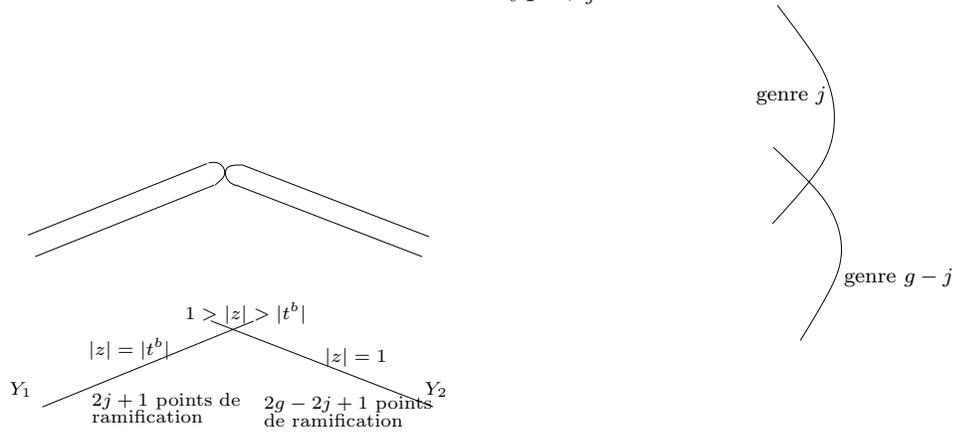}}
\end{figure}

\begin{figure}
\caption{Point de type $\alpha_j$}
\label{PointAlphaj}
\centerline{\input{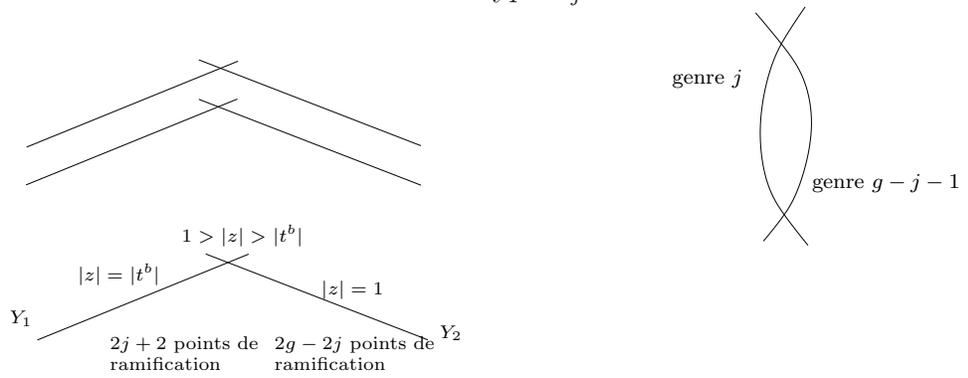}}
\end{figure}

 On peut lire les normes $|c_i-c_k|$ en fonction des composantes sur lesquelles se sp\'ecialisent $b_i$ et $b_k$ dans $Y_s$. On trouve directement 

\centerline{\begin{tabular}{ll}
type $\alpha_j$ & $|\Delta|=|A|^{4g+2}\prod_{i \not = k}|c_i-c_k|=|A|^{4g+2}|t^{b (2j+1) (2j+2)}|$ \\
type $\beta_j$ & $|\Delta|=|A|^{4g+2}\prod_{i \not = k}|c_i-c_k|=|A|^{4g+2}|t^{b 2j (2j+1)}|$ 
\end{tabular}}

On construit alors une base de $\H^0(\C, \omega_{\C/\Spec(R)})$ comme dans \cite{Kau}, proposition 5.5.
\begin{Lem}
Dans le cas $\beta_j$ on a $$div_{vert}({dx \over 2y+Q})=-{1 \over 2}\nu(A) Y_1+\left ( -{1\over 2}\nu(A)-{1 \over 2} (2j-1)b\right )Y_2$$

et dans le cas $\alpha_j$  on a
$$div_{vert}({dx \over 2y+Q})=-{1 \over 2}\nu(A) Y_1+\left ( -{1\over 2}\nu(A)-{1 \over 2} (2j)b\right )Y_2$$
\end{Lem}

\noindent{\it Preuve :}  On donne la d\'emonstration dans le cas $\beta_j$, la d\'emonstration \'etant semblable dans le cas $\alpha_j$.
Soit $i\in \{1, 2\}$ fix\'e. Soit $x'$ un param\`etre de $Y_\eta$ relevant un param\`etre de $Y_i$ (on peut prendre par exemple $x$ pour $Y_1$ et $xt^{-b}$ pour $Y_2$). Notons $U=\Spec(A)$ un voisinage affine dans $Y$ du point g\'en\'erique de $Y_i$ ne rencontrant pas l'autre composante. Notons $B$ la cl\^oture int\'egrale de $A$ dans $\kappa(\C)$. D'apr\`es \cite{Liu}, Lemme 1, $B$ est libre sur $A$, de rang 2 et il existe $y'$ tel que $\{1, y'\}$ soit une base de $B$ sur $A$. Notons ${y'}^2+\tilde Q(x')y'-\tilde P(x')$ une \'equation affine associ\'ee. Notons
$$\omega={dx'\over (2y'+\tilde Q(x'))}$$
 une base locale du faisceau dualisant (c.f. \cite{Liu2}, Lemme 2). 

D'autre part, comme $\{1, y\}$ est une base de $A \times K$ on a $y=t^{\zeta_i} y'+s(x)$. Notons $\tilde f=\tilde Q^2+4\tilde P$. En rempla\c cant $y$ par $t^{\zeta_i} y'+s(x)$ dans  l'\'equation $y^2+Qy-P$ on trouve que 
$$t^{\zeta_i} \tilde Q(x')=2s(x)+Q(x)$$ 
$$-t^{2\zeta_i} \tilde P(x')=Q(x)s(x)+s(x)^2-P(x)$$
donc $t^{2\zeta_i}\tilde f(x')=f(x)$. De plus, comme $Y_i$ est r\'eduite (c'est \`a dire que la r\'eduction de $\tilde f$ \`a $k$ est non nulle), on a $|\tilde f|=1$

Notons $\nu_{Y_1}(f)=\sup_{|z|=1}(\nu(f(z)))$ et $\nu_{Y_2}(f)=\sup_{|z|=|t^b|}(\nu(f(z)))$. On a alors  $\nu_{Y_1}(f)=\nu(A)$ et  $\nu_{Y_2}(f)=\nu(A)+(2j+1)b$.

On a donc $2\zeta_1=\nu(A)$ et $2\zeta_2=\nu(A)+(2j+1)b$. 

Finalement, comme $${dx \over 2y+Q} = {t^b dx' \over {t^{\zeta_i}} (2y'+\tilde Q(x'))}$$ et $\omega$ est une base locale, on trouve que 
$$\nu_{Y_1}({dx \over 2y+Q})=-{1 \over 2} \nu(A)$$ et $$\nu_{Y_2}({dx \over 2y+Q})=b-{1 \over 2} \nu(A)-{1\over 2}(2j+1)b.$$ 
Ce qui nous donne bien le r\'esultat souhait\'e. \Eproof
\medskip

Notons $\beta$ l'un des points $b_i$ qui se sp\'ecialise sur la composante $Y_2$ et, pour $0 \le i \le g-1$, posons 
$$\omega_i=(x-\beta)^i {dx \over 2y+Q}.$$
Par suite, d\'efinissons $m_i=min(\nu_{Y_2}(\omega_i), \nu_{Y_1}(\omega_i))$. On a alors, dans le cas pr\'esent $$m_i=min(\nu_{Y_2}(\omega_i), \nu_{Y_1}(\omega_0))=min(ib+\nu_{Y_2}(\omega_0), \nu_{Y_1}(\omega_0)).$$

\begin{Lem}
Les $\{t^{-m_i} \omega_i\}$ forment une base de $\H^0(\C, \omega_{\C/\Spec(R)})$
\end{Lem}

\noindent{\it Preuve :}  Comme les $\omega_i$ forment une base de $\H^0(\C_\eta, \omega_{\C_\eta / K})$, ils sont lin\'eai\-rement ind\'ependant. Le diviseur $div(t^{-m_i} \omega_i)$ est effectif par construction donc $t^{-m_i} \omega_i \in \H^0(\C, \omega_{\C/\Spec(R)})$. Il reste \`a montrer que les $\omega_i$ sont g\'en\'erateurs. 

Soit $\omega \in \H^0(\C, \omega_{\C/\Spec(R)})$, $\omega$ est de la forme $t^e h(x) {dx \over 2y+Q}$ avec $h \in R[X]$ primitif de degr\'e $d \le g-1$.

Comme $div(\omega) \ge 0$, on a pour tout $k \in \{1, 2\}$ $$\nu_{Y_k}(\omega)=e+\nu_{Y_k}(h(x))+\nu_{Y_k}({dx \over 2y+Q}) \ge 0.$$
Or, comme $h$ est primitif, $\nu_{Y_2}(h) \le db$ et $\nu_{Y_1}(h) \ge 0$ donc $e+db+\nu_{Y_2}(\omega_0) \ge 0$ et $e +\nu_{Y_1}(\omega_0)\ge 0$. C'est \`a dire que $e\ge -m_d$. Par suite, il existe un \'el\'ement $c\in R$ tel que $\omega-ct^{-m_d}\omega_d=h_1(x) {dx \over 2y+Q}$ soit dans $\H^0(\C, \omega_{\C/\Spec(R)})$ avec $h_1$ de degr\'e inf\'erieur ou \'egal \`a $d-1$. Par r\'ecurrence, on montre ainsi que la famille $\{t^{-m_i}\omega_i\}$ est g\'en\'eratrice. \Eproof
\medskip

Finalement, on trouve que 
$$\Ord_s({dx \over 2y+Q} \wedge \ldots \wedge {x^{g-1}dx \over 2y+Q})=\sum m_i$$

Un calcul simple nous montre que dans le cas pr\'esent on a 

\centerline{\begin{tabular}{ll}
type $\alpha_j \Rightarrow$ & $\Ord_s \left({dx \over 2y+Q} \wedge \ldots \wedge {x^{g-1}dx \over 2y+Q}\right)=-{j^2 +j \over 2}b$ \\
type $\beta_j \Rightarrow$ & $\Ord_s \left({dx \over 2y+Q} \wedge \ldots \wedge {x^{g-1}dx \over 2y+Q}\right)=-{j^2 \over 2}b$
\end{tabular}}

Finalement, on trouve que 

\centerline{\begin{tabular}{ll}
type $\alpha_j \Rightarrow $ & $\Ord_s(\Lambda)= 2(g-j)(j+1)b$ \\
type $\beta_j \Rightarrow $ & $\Ord_s(\Lambda)=2(g-j)jb$
\end{tabular}}

Il ne reste donc plus qu'\`a relier $\delta_s$ \`a $b$.

Dans le cas $\beta_j$, on a $2 \delta_s=b$, ceci provient du morphisme
$$\Spec R[[xy]]/(xy-t^b) \rightarrow \Spec R[[u, v]]/(uv-t^{2b})$$ donn\'e par $u\mapsto x^2$ et $v \mapsto y^2$.

Dans le cas $\alpha_j$, on trouve directement que $\delta_s = 2b$ (car les points singuliers de $\C_s$ sont \'echang\'es par $\sigma$). Finalement, on trouve le r\'esultat annonc\'e. \Eproof

\section{Calcul en caract\'eristique 2}

Dans cet appendice, on se propose de montrer les in\'egalit\'es (\ref{TheoPrinc}) 
par des calculs explicites dans le cas de certaines courbes en caract\'eristique r\'esiduelle 2. Les r\'esultats ne sont pas utilis\'es dans le reste de l'article.

\begin{Pro}
Soit $R$ un anneau de valuation discr\`ete de caract\'eristiques $(0, 2)$.
Soit $\C \rightarrow \Spec(R)$ une courbe stable hyperelliptique dont la fibre g\'en\'erique est lisse et dont la fibre sp\'eciale contient au plus deux composantes irr\'eductibles. On a $g^2 \delta_s \ge \Ord_s(\Lambda) \ge g \delta_s$.
\end{Pro}

\noindent{\it Preuve :} Le calcul est essentiellement le m\^eme que dans le cas de la caract\'eristique r\'esiduelle $\not = 2$ sauf pour le calcul du discriminant, on se contente donc de d\'etailler celui-ci.

Si $\C$ est lisse, on a $\Ord_s(\Lambda)=\delta_s=0$ (voir par exemple \cite{Kau}, proposition 6.3). Si la fibre sp\'eciale de $\C$ poss\`ede une unique composante irr\'eductible, celle-ci poss\'edant un point singulier de type $\alpha_0$, on se ram\`ene au cas de deux composantes par \'eclatement (apr\`es extension eventuelle de $R$). 

Notons $\sigma$ l'involution hyperelliptique. Soit $y^2+Q(x)y-P(x)$ une \'equation affine de $\C_\eta$. Le discriminant de l'\'equation est alors $f=Q^2+4P$. La caract\'eristique r\'esiduelle \'etant 2, on voit que les points fixes de $\C_\eta$ sous l'action de $\sigma$ se sp\'ecialisent par paire (exactement par paire pour des raisons semblables \`a celles expos\'ees dans la d\'emonstration du lemme \ref{SpecialisationSimple}).

Notons $\C'$ la r\'esolution compl\`ete des sp\'ecialisations des points fixes de $\C_\eta$ sous l'action de $\sigma$. On peut alors avoir pr\'ecis\'ement la r\'epartition des sp\'ecialisations des points fixes grace \`a \cite{Hen_p_adique}, th\'eor\`eme 4.5, en fonction du conducteur et de l'\'epaisseur des points doubles dans $\C/\gsigma$. On donne aux figures \ref{Fig1}, \ref{Fig2}, \ref{Fig3} et \ref{Fig4} les configurations possibles  pour $(\C'/\gsigma)_s$ dans le cas ou les composantes sont de genres $j$ et $g-j$ (respectivement $j$ et $g-j-1$ dans le cas d'une courbe ayant deux points singuliers de type $\alpha_j$). On a repr\'esent\'e sur chaque figure la fibre sp\'eciale de $\C$ ainsi que celle de $\tilde\C$ (o\`u $\tilde\C$ est la courbe obtenue de $\C$ en d\'eployant les sp\'ecialisations des points fixes de $\C$ sous l'action de $\sigma$). Notons $b$ l'\'epaisseur du point singulier de $\C'/\gsigma$.

\begin{figure}
\caption{La r\'eduction poss\`ede un seul point singulier et $b > 2\nu_K(2)$}
\centerline{\input{type1.pstex_t}}
\label{Fig1}
\end{figure}

Dans le cas des figures \ref{Fig1}, \ref{Fig2} et \ref{Fig3} (correspondant aux r\'eduction possibles des courbes stables hyperelliptiques dont la fibre sp\'eciale poss\`ede deux composantes irr\'eductibles et ayant un unique point singulier, celui-ci \'etant de type $\beta_j$), on trouve, par un calcul direct,  $|disc(f)|=|4^{2g+2} t^{b 2j(2j+1)}|$ et par suite, d'apr\`es les calculs fait pr\'ecedement pour les cas $\beta_j$ (c.f. calcul en caract\'eristique z\'ero), on a $$\Ord_s(\Lambda)=2(g-j)jb=(g-j)j\delta_s.$$

\begin{figure}
\caption{La r\'eduction poss\`ede un seul point singulier et $b=2\nu_K(2)$}
\centerline{\input{type2.pstex_t}}
\label{Fig2}
\end{figure}

\begin{figure}
\caption{La courbe poss\`ede un unique point singulier et $b < 2\nu_K(2)$}
\centerline{\input{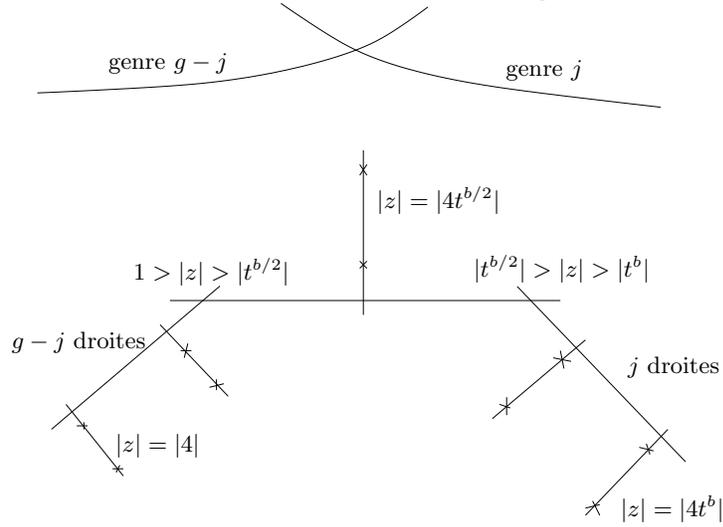}}
\label{Fig3}
\end{figure}
\begin{figure}
\caption{La courbe poss\`ede deux points singuliers de type $\alpha_j$}
\centerline{\input{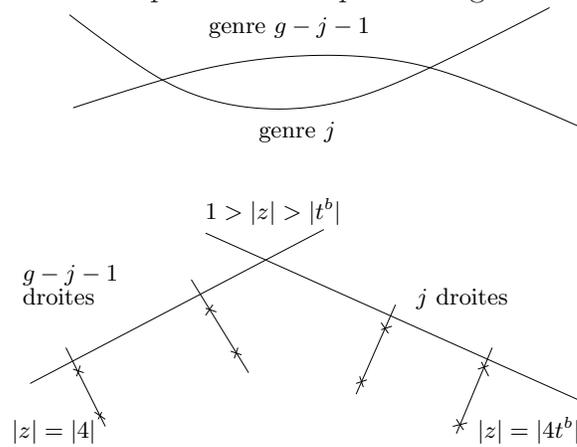}}
\label{Fig4}
\end{figure}

Dans le cas de la figure \ref{Fig4} (correspondant au courbes stables hyperelliptiques dont la fibre sp\'eciale poss\`ede deux points singuliers de type $\alpha_j$). Des calculs analogues au cas de la caract\'eristique z\'ero montrent que  $$|disc(f)|=|4^{2g+2} t^{b (2j+2)(2j+1)}|,$$ on a donc $$\Ord_s(\Lambda)=2(g-j)(j+1)b+b=((g-j)j)\delta_s$$

\clearpage
\bibliographystyle{amsplain}

\bibliography{ModHypLiu.bbl}

\begin{flushleft}
Sylvain Maugeais \\
Laboratoire A2X, UMR 5465 CNRS \\
Universit\'e de Bordeaux I \\
351 cours de la Lib\'eration, \\
33405 Talence Cedex, France \\
e-mail : {\tt maugeais@math.u-bordeaux.fr}
\end{flushleft}
\end{document}